\documentclass{amsart}
\usepackage{amssymb}
\usepackage[all,cmtip]{xy}
\DeclareMathOperator{\Val}{Val}
\DeclareMathOperator{\Gr}{Gr}
\DeclareMathOperator{\AGr}{\overline{Gr}}
\DeclareMathOperator{\vol}{vol}
\DeclareMathOperator{\Kl}{Kl}
\DeclareMathOperator{\Flag}{Flag}
\DeclareMathOperator{\Hom}{Hom}
\DeclareMathOperator{\Dens}{Dens}
\newtheorem{Proposition}{Proposition}[section]
\newtheorem{Theorem}[Proposition]{Theorem}
\newtheorem{Corollary}[Proposition]{Corollary}
\newtheorem{Definition}[Proposition]{Definition}
\newcommand\flag[2]{\left[\begin{array}{c} #1\\ #2 \end{array}
  \right]} 

\title{Algebraic integral geometry}
\author{Andreas Bernig}

\email{bernig@math.uni-frankfurt.de}

\address{Institut f\"ur Mathematik, Goethe-Universit\"at Frankfurt,
Robert-Mayer-Str. 10, 60054 Frankfurt, Germany}

\begin{document}

\begin{abstract}
A survey on recent developments in (algebraic) integral geometry is given. The
main focus lies on algebraic structures on the space of translation invariant
valuations and applications in integral geometry.
\end{abstract}

\maketitle 

\tableofcontents

\section{Algebraic integral geometry}

Algebraic integral geometry is a relatively modern part of integral geometry. It aims at proving geometric formulas (kinematic formulas, Crofton formulas, Brunn-Minkowski-type inequalities etc.) by taking a structural viewpoint and employing various algebraic techniques, including abstract algebra, Lie algebras and groups,  finite- and infinite-dimensional representations, classical invariant theory, Gr\"obner bases, cohomology theories, algebraic geometry and so on.  

The situation can be roughly compared to symbolic integration. In order to integrate a given (say sufficiently elementary) function, there is no need to know anything about the definition of the integral. It suffices to know a certain number of integration rules, like partial integration and the substitution rule. In algebraic integral geometry, the corresponding rules for computing geometric integrals are worked out. The fundamental theorem of algebraic integral geometry is one of these rules. 

The main object of the theory is the {\it space of all translation invariant valuations}. Here the emphasis is on {\it space}, since in general not a single valuation but the set of all valuations is studied. Roughly speaking, this space is a graded commutative algebra satisfying Poincar\'e duality and Hard Lefschetz theorem. Moreover, there is a Fourier transform and a convolution product on this algebra and all these algebraic structures reflect geometric properties and formulas. 

Among the most important contributions to algebraic integral geometry are Nijenhuis observation \eqref{eq_kf_nijenhuis} (who moreover speculated about a possible algebraic structure explaining it), the theorems by Hadwiger (Theorem \ref{thm_hadwiger}), P. McMullen (Theorem \ref{thm_mcmullen}), Klain (Subsection \ref{subsec_klain_embedding}) and Schneider (see Subsection \ref{subsec_schneider_embedding}). A spectacular breakthrough was achieved by Alesker in 2001, who proved the {\it McMullen conjecture} (in fact a much stronger version of it, see Subsection \ref{subsec_irr_thm}) and subsequently introduced many of the algebraic structures mentioned above. 

The structure of the present paper is as follows. 

After a short reminder of some classical integral-geometric formulas in Section \ref{sec_classical_intgeo}, we will explain the new algebraic tools in Section \ref{sec_alg_structures}. 

The transition between algebra and geometry is done in Section \ref{sec_ftaig}, where the theoretical background for integral geometry of subgroups of $\mathrm{SO}(n)$ is given. 

In Section \ref{sec_un}, this program is carried out in a special and  important case, namely for the group $G=\mathrm{U}(n)$, yielding {\it hermitian integral geometry}. Section \ref{sec_other_groups} gives an overview of integral geometry for other groups and three important open problems are stated in Section \ref{sec_open_problems}. 

The reader is invited to read J.~Fu's survey \cite{fu_aig} which has some non-empty intersection with the present paper. 

\subsection*{Acknowledgments} I was happy to profit from many discussions with and talks by Semyon Alesker and Joseph Fu on algebraic integral geometry. The present text is strongly influenced by their ideas and I am grateful to them. The terms {\it Algebraic integral geometry} and {\it Fundamental theorem of algebraic integral geometry} were invented by Fu. I also thank Gautier Berck, Franz Schuster and Christoph Th\"ale for numerous useful remarks on this text.  

\section{Classical integral-geometric formulas}
\label{sec_classical_intgeo}

Let us fix some notations for the rest of the paper. The $n$-dimensional unit ball is denoted by $B$. Let $\omega_n$ be its volume. The {\it flag coefficients} are defined by 
\begin{displaymath}
 \flag{n}{k}:=\binom{n}{k}\frac{\omega_n}{\omega_k \omega_{n-k}}.
\end{displaymath}

For an odd number $2m+1$, we set 
\begin{displaymath}
 (2m+1)!!=1 \cdot 3 \cdot 5 \cdot \ldots \cdot (2m+1)
\end{displaymath}
and use the convention $(-1)!!=1$. 

If $V$ is a finite-dimensional Euclidean vector space and $G$ is any subgroup of $\mathrm{SO}(V)$, we let $\bar G$ be the group generated by $G$ and translations. This group has a canonical measure, which is the product
of the Haar probability measure on $G$ and the Lebesgue measure on the
translation group. In particular, the measure of the set $\{\bar g \in \bar G:
\bar g(x) \in K\}$, $x \in V$ equals the volume of $K$ for every compact convex set $K$. 

We will use the following terminology: {\it subspaces} will always be {\it linear subspaces}, while {\it planes} will always be {\it affine planes}. 
The {\it Grassmann manifold} of all $k$-dimensional subspaces in $V$ is denoted by $\Gr_kV$. The {\it affine Grassmann manifold} of all $k$-planes  is denoted by $\AGr_kV$. 

\subsection{Valuations}

Throughout this paper, $V$ denotes a finite-dimensional vector space. The
space of non-empty compact convex subsets in $V$ is denoted by $\mathcal{K}(V)$. With respect to {\it Minkowski addition} 
\begin{displaymath}
 K+L=\left\{ x+y| x \in K, y \in L \right\},
\end{displaymath}
$\mathcal{K}(V)$ is a semigroup. This space has a natural topology, called the {\it Hausdorff-topology} which is
defined as follows: 
\begin{displaymath}
 d_H(K,L):=\inf_{r \geq 0} \left\{K \subset L+rB, L \subset K+rB\right\},
\quad K,L \in \mathcal{K}(V). 
\end{displaymath}
Here $B$ is the unit ball for some euclidean scalar product. The metric $d_H$ depends on the choice of this scalar product, but the induced
topology does not.

\begin{Definition} \label{def_valuation}
Let $A$ be a semigroup. A functional $\mu:\mathcal{K}(V) \to A$ is called a
valuation if 
\begin{displaymath}
 \mu(K \cup L)+\mu(K \cap L)=\mu(K)+\mu(L)
\end{displaymath}
whenever $K,L, K \cup L \in \mathcal{K}(V)$. 
\end{Definition}

The case of $A=\mathbb{R}$ (or $A=\mathbb{C}$) is the most important one. Everyone is familiar with
at least two examples of real-valued valuations. The first one is the constant
valuation $\mu(K)=1$ for all $K \in \mathcal{K}(V)$. This valuation is called
{\it Euler characteristic} and denoted by $\chi$. The name needs some
explanation: in fact there is a canonical way to extend this valuation to finite
unions of compact convex sets, and this extension equals the Euler
characteristic with respect to Borel-Moore homology.  

The second familiar example of a valuation is the volume, which we denote by
$\vol$. Note that this valuation depends on the choice of a Euclidean
metric on $V$ (or at least on the choice of a Lebesgue measure).  

A particularly important class of valuations is that of continuous, translation
invariant valuations. The valuation $\mu$ is called translation invariant if
$\mu(K+t)=\mu(K)$ for all $t \in V$. Euler characteristic and volume clearly
have this property. We will see later on that all continuous, translation
invariant valuations arise in some way from these two basic ones.  

\begin{Definition}
The space of complex-valued, continuous, translation invariant valuations is denoted by $ \Val$. If $G$ is a subgroup of $\mathrm{GL}(V)$, then 
\begin{displaymath}
 \Val^G(V):=\{\mu \in \Val| \mu(gK)=\mu(K) \quad \forall g \in G\}.
\end{displaymath}
\end{Definition}

Before studying the space $\Val$, let us give some
other important examples of valuations, which do not belong to $\Val$. 

First of all, non-continuous valuations (on the space of polytopes), the {\it Dehn
functionals}, played a
central role in Dehn's solution of Hilbert's 3rd problem. Another famous non-continuous example is the {\it affine surface area}, which is semi-continuous (see \cite{lutwak91} and \cite{ludwig_reitzner99} and the references therein for more information). 

Since $\mathcal{K}(V)$, endowed with the Minkowski addition, is a semigroup,  
we may take $A=\mathcal{K}(V)$ in Definition
\ref{def_valuation}. It is easy to see that $\mu(K)=K$ defines a valuation.
More interesting examples are the {\it intersection body operator} (defined on a subset of $\mathcal{K}(V)$) and the {\it projection body operator}. See \cite{ludwig06, ludwig06_survey, schuster08, abardia_bernig} for more information.  

If $A=\mathop{Sym}^*V$, the space of symmetric tensors over $V$, then $A$-valued valuations are called {\it tensor valuations}. Their study
has been initiated by McMullen \cite{mcmullen97} and Alesker \cite{alesker99}. Recently, 
remarkable progress in the study of kinematic formulas for tensor valuations
was made by Hug, Schneider and R.~Schuster \cite{hug_schneider_schuster_a, hug_schneider_schuster_b}. One may hope and expect that some algebraic tools will be useful in simplifying their formulas.   
  
\subsection{Intrinsic volumes}
\label{subsec_int_vol}

Let $V$ be a Euclidean vector space of dimension $n$. At the heart of integral geometry are the {\it intrinsic volumes} $\mu_0,\ldots,\mu_n$. We give four equivalent definitions. 

First of all, we may use {\it projections} onto lower-dimensional subspaces. For $0 \leq k \leq n$, the group $\mathrm{SO}(V)$ acts transitively on the Grassmannian $\Gr_k(V)$ of $k$-dimensional subspaces in $V$. We endow this manifold with the unique invariant probability measure $dL$. Then 
\begin{equation} \label{eq_int_vol_projections}
 \mu_k(K):=\flag{n}{k} \int_{\Gr_k(V)} \vol_k(\pi_LK)dL
\end{equation}
defines an element $\mu_k \in \Val^{\mathrm{SO}(V)}$. Here $\vol_k$ denotes the $k$-dimensional Lebesgue measure on the subspace $L$ and $\pi_LK$ is the orthogonal projection of $K$ onto $L$. This formula (and some more general versions) is called {\it Kubota formula}. 

For the second definition, we use {\it intersections} instead of projections. We let $\AGr_k(V)$ be the $k$-dimensional affine Grassmannian on which we use the unique  $\overline{\mathrm{SO}(V)}$-invariant measure $dE$
such that the measure of planes intersecting the unit ball equals $\omega_{n-k}$. Then we set 
\begin{equation} \label{eq_int_vol_intersections}
 \mu_k(K):=\flag{n}{k} \int_{\AGr_{n-k}(V)} \chi(K \cap E) dE.
\end{equation}
The equivalence of this definition with the previous one is an elementary exercise. Formula \eqref{eq_int_vol_intersections} is called {\it Crofton formula}. More general Crofton formulas play an important role in algebraic integral geometry, see Subsection \ref{subsec_product}.  

The third description of the $\mu_k$ is rather a characterization than a definition. Looking at the $\mu_k$ defined as above, one sees that 
\begin{enumerate}
 \item $\mu_k$ is a continuous, translation invariant and $\mathrm{SO}(V)$-invariant
valuation,
\item $\mu_k$ is of degree $k$, i.e. $\mu_k(tK)=t^k\mu_k(K)$ for all $t \geq 0$ and
\item the restriction of $\mu_k$ to a $k$-plane equals the $k$-dimensional Lebesgue measure on that plane. 
\item $\mu_k$ is even, i.e. $\mu_k(-K)=\mu_k(K)$.
\end{enumerate}

In fact, the $\mu_k$ are uniquely characterized by these properties, as we will see in Subsection \ref{subsec_klain_embedding}.

Some of the $\mu_k$ are well-known: $\mu_0$ is the Euler characteristic $\chi$ which was mentioned in the introduction. $\mu_n$ is the usual Lebesgue measure, $\mu_{n-1}$ is half of the surface area and $\mu_1$ is a constant times the {\it mean width}. 

The intrinsic volumes may also be defined by the {\it Steiner formula}. For $t \geq 0$, let $K+tB$ be the $t$-tube around $K$. Then $\vol(K+tB)$ turns out to be a polynomial in $t$ given by
\begin{equation} \label{eq_steiner}
 \vol(K+tB)=\sum_{k=0}^n \mu_{n-k}(K) \omega_k t^k.
\end{equation}

Taking $K=B$, we easily get 
\begin{equation}
 \mu_k(B)=\binom{n}{k}\frac{\omega_n}{\omega_{n-k}}.
\end{equation}

\subsection{Kinematic formulas}
\label{subsec_kinformula_son}

The most important formulas of integral geometry are the {\it kinematic
formulas}:
\begin{equation} \label{eq_kin_formula_son}
 \int_{\overline{\mathrm{SO}(V)}} \mu_i(K \cap \bar g L)d\bar
g=\sum_{k,l} c_{k,l}^i \mu_k(K)\mu_l(L), 
\end{equation}
where 
\begin{displaymath}
c_{k,l}^i=\left\{\begin{array}{c l}\flag{n+i}{i} \flag{n+i}{k}^{-1} & k+l=n+i\\ 0 & k+l \neq n+i. \end{array}\right.
\end{displaymath}

It may be checked that the constants on the right hand side are correct by
plugging in balls of different radii (template method, see below). 

Looking at the formula, one makes the following observations. Since 
\begin{displaymath}
 \flag{n+i}{k}=\flag{n+i}{n+i-k},
\end{displaymath}
the coefficients on the right hand side are
symmetric, i.e. $c_{k,l}^i=c_{l,k}^i$. This reflects of course the fact that changing the role of $K$ and
$L$ in the integral on the left hand side of the formula does not change its
value. Next, we observe that the {\it total degree} $n+i$ at the right hand side
is the degree $i$ on the left hand side plus the dimension of the ambient
space. Another symmetry property for the coefficients comes from Fubini's theorem:  
\begin{displaymath}
 \int_{\overline{\mathrm{SO}(V)}} \int_{\overline{\mathrm{SO}(V)}} \mu_i(K \cap \bar g L \cap
\bar h M) d\bar g d\bar h=\int_{\overline{\mathrm{SO}(V)}} \int_{\overline{\mathrm{SO}(V)}}
\mu_i(K \cap \bar g L \cap
\bar h M) d\bar h d\bar g.
\end{displaymath}
This translates to 
\begin{displaymath}
 \sum_r c_{r,m}^ic_{k,l}^r=\sum_r c_{r,l}^ic_{k,m}^r.
\end{displaymath}

Nijenhuis \cite{nijenhuis74}
made a less obvious observation: 
Renormalizing 
\begin{displaymath}
 \tilde \mu_k:=\frac{\pi^n k!\omega_k}{\pi^k n!\omega_n}\mu_k,
\end{displaymath}
the
kinematic formula \eqref{eq_kin_formula_son} becomes 
\begin{equation} \label{eq_kf_nijenhuis}
 \int_{\overline{\mathrm{SO}(V)}} \tilde \mu_i(K \cap \bar g L)d\bar g=\sum_{k+l=n+i}
\tilde \mu_k(K)\tilde \mu_l(L).
\end{equation}
Hence all coefficients on the right hand side become $1$.  

At first glance, this may seem to be trivial, since we may change the constants
on the right hand side to whatever we want by rescaling the $\mu_k$. However, a
closer look reveals that we only have $n+1$ free parameters (one for the
scaling of each
$\mu_k$), but $\binom{n+2}{2}$ non-zero coefficients on the right hand side. Nijenhuis
speculated that there exists some algebraic structure explaining this strange
fact (``...the suggestion of an underlying algebra with the $c$'s as structure constants was inevitable'' \cite{nijenhuis74}). It turns out that this is indeed the case, as we will see below.

An array of {\it additive kinematic formulas} arises if we replace intersection
by Minkowski addition:
\begin{equation} \label{eq_add_kf}
 \int_{\mathrm{SO}(V)} \mu_i(K +g L)dg=\flag{2n-i}{n-i}\sum_{k+l=i}
\flag{2n-i}{n-k}^{-1}\mu_k(K)\mu_l(L).
\end{equation}
In this case, there is an analogous statement as in Nijenhuis' observation:
after renormalizing 
\begin{displaymath}
 \tilde \mu_k:=\frac{(n-k)!\omega_{n-k}}{n!\omega_n} \mu_k, k=0,\ldots,n,
\end{displaymath}
the additive kinematic formula \eqref{eq_add_kf} reads 
\begin{displaymath}
 \int_{\mathrm{SO}(V)} \tilde \mu_i(K +g L)dg=\sum_{k+l=i}
\tilde\mu_k(K)\tilde\mu_l(L). 
\end{displaymath}
We will see later an explanation of this fact too. It will also turn out that
the usual and the additive kinematic formula are {\it dual to
each other} in a precise sense and one can be derived from the other. 

\subsection{Hadwiger's theorem}
\label{subsec_hadwiger}

We have already seen that $\mu_k \in \Val^{\mathrm{SO}(V)}$. Hadwiger's theorem states
conversely that {\it all} valuations in $\Val^{\mathrm{SO}(V)}$ are obtained by linear
combinations of intrinsic volumes.

\begin{Theorem} \label{thm_hadwiger}
 The vector space $\Val^{\mathrm{SO}(V)}$ of continuous, translation invariant,
$\mathrm{SO}(V)$-invariant valuations on a Euclidean vector space $V$ of dimension $n$
has the basis
\begin{displaymath}
 \mu_0,\mu_1,\ldots,\mu_n.
\end{displaymath}
\end{Theorem}

An elementary and nice proof may be found in \cite{klain_rota}. 

Hadwiger's theorem is quite powerful. It leads to all formulas which we have
stated before. Indeed, let us look for instance at the additive kinematic
formula \eqref{eq_add_kf}. For each fixed body $L$, the left hand side of this formula is a
valuation in $K$. It is easy to prove that this valuation belongs to
$\Val^{\mathrm{SO}(V)}$, hence it may be written in the form $\sum_{k=0}^n
d_k(L)\mu_k(K)$. Next, fixing $K$, one easily gets that $d_k$ is also an
element of $\Val^{\mathrm{SO}(V)}$ for each fixed $k$, hence $d_k(L)=\sum_{l=0}^n
d_{kl}\mu_l(L)$ with complex numbers $d_{kl}$. We thus know that 
\begin{displaymath}
 \int_{\mathrm{SO}(V)} \mu_i(K +g L)dg=\sum_{k,l=0}^n d_{kl}^i
\mu_k(K)\mu_l(L)
\end{displaymath}
for some fixed constants $d_{kl}^i$. There is a nice trick to determine these
constants, which is called the {\it template method}. We plug in on both sides of
the equation special convex bodies $K$ and $L$ for which we may compute
the integral on the left hand side and the intrinsic volumes on the right
hand side to obtain a system of linear
equations on the $d_{kl}^i$. Solving this system yields the $d_{kl}^i$. More precisely, let us take $K=rB, L=sB$ (where $B$ is as always the unit ball). The left hand side equals $\mu_i((r+s)B)=(r+s)^i\binom{n}{i}\frac{\omega_n}{\omega_{n-i}}$. The right hand side equals 
\begin{displaymath}
 \sum_{k,l} d_{k,l}^i r^k \binom{n}{k}\frac{\omega_n}{\omega_{n-k}} s^l \binom{n}{l}\frac{\omega_n}{\omega_{n-l}}. 
\end{displaymath}

Comparing the coefficients of $r^j s^{i-j}$ on both sides gives us 
\begin{displaymath}
 \binom{i}{j} \binom{n}{i}\frac{\omega_n}{\omega_{n-i}}=d_{j,i-j}^i \binom{n}{j}\frac{\omega_n}{\omega_{n-j}} \binom{n}{i-j}\frac{\omega_n}{\omega_{n-i+j}},
\end{displaymath}
which is \eqref{eq_add_kf}.

\subsection{General Hadwiger theorem} 

The last theorem from classical integral geometry which we want to mention is
the {\it general Hadwiger theorem}. It applies to more general valuations, but
we will state (and prove) it only in the special case of translation
invariant valuations. 

\begin{Theorem} \label{thm_gen_hadwiger}
 Let $\phi \in \Val$. Then 
\begin{displaymath}
 \int_{\overline{\mathrm{SO}(V)}} \phi(K \cap \bar g L) d\bar g=\sum_{l=0}^n c_l(K)
\mu_l(L),
\end{displaymath}
where 
\begin{displaymath}
 c_l(K):=\int_{\AGr_{n-l}} \phi(K \cap E) dE.
\end{displaymath}
\end{Theorem}

The theorem can be proved using Hadwiger's characterization theorem and a limit
argument. We will give another, more conceptual proof, in Subsection \ref{subsec_ftaig}.

\section{Algebraic structures on valuations}
\label{sec_alg_structures}

The main object of algebraic integral geometry is the space $\Val=\Val(V)$ of
continuous, translation invariant valuations on an $n$-dimensional vector space
$V$. This space has a surprisingly rich algebraic structure which we are going
to describe in this section. 

\subsection{McMullen's decomposition}

A valuation $\mu$ is {\it of degree $k$} if $\mu(tK)=t^k\mu(K)$ for all
$t \geq 0$ and all $K$. It is {\it even} if $\mu(-K)=\mu(K)$ and {\it
odd} if $\mu(-K)=-\mu(K)$. The corresponding subspaces of $\Val$ are denoted by $\Val_k^+, \Val_k^-$. 

McMullen \cite{mcmullen77} proved the following decomposition: 
\begin{Theorem} \label{thm_mcmullen}
 \begin{equation} \label{eq_mcmullen}
  \Val=\bigoplus_{\substack{k=0,\ldots,n\\ \epsilon = \pm}} \Val_k^\epsilon. 
 \end{equation}
\end{Theorem}

In particular, the degree of a valuation is an integer between $0$ and the
dimension of the ambient space. We refer to \eqref{eq_mcmullen} as the {\it
McMullen grading}. 

McMullen's theorem allows us to introduce a Banach space structure on $\Val$ by setting
\begin{displaymath}
 \|\mu\|:=\sup \left\{|\mu(K)|: K \in \mathcal{K}(V), K \subset B \right\}
\end{displaymath}
Then $(\Val,\|\cdot\|)$ is a Banach space. Choosing another scalar product on $V$ gives an equivalent norm. Hence we get a uniquely defined {\it Banach space structure} on $\Val$. 

\subsection{Klain embedding}
\label{subsec_klain_embedding}

We now suppose that we have a fixed Euclidean structure on $V$. This is not
strictly necessary but simplifies the exposition. 

Klain found a nice way to describe {\it even} continuous, translation invariant
valuations. In \cite{klain00}, he first characterized the volume as the only valuation (up to a
multiplicative constant) which is continuous, translation invariant, even and
{\it simple} (i.e. vanishing on lower-dimensional sets). 

Klain's characterization of the volume implies that given $\mu \in \Val_k^+$
and a $k$-dimensional subspace $E$, the restriction $\mu|_E$ is a multiple
$\Kl_\mu(E)$ of the $k$-dimensional volume on $E$. Indeed, this follows once we know that $\mu|_E$ is simple. If $F \subset E$ is a subspace of minimal dimension such that $\mu|_F \neq 0$, then $\mu|_F$ is simple and hence a multiple of the volume on $F$. Since $\mu$ is of degree $k$, this is only possible if $E=F$.  

The continuous function 
\begin{equation} \label{eq_klain_embedding}
 \Kl_\mu:\Gr_k(V) \to \mathbb{C}
\end{equation}
is called the {\it Klain function} of $\mu$. The induced map
\begin{displaymath}
\Kl: \Val_k^+  \hookrightarrow C(\Gr_kV)
\end{displaymath}
is called the {\it Klain embedding}. To see that this map is indeed injective, we
suppose that $\Kl_\mu=0$ for some $\mu \in \Val_k^+$. Then the restriction of
$\mu$ to any $(k+1)$-dimensional subspace $F$ is simple, hence a multiple of the
$(k+1)$-dimensional volume on $F$. Since $\mu$ is $k$-homogeneous, this is only
possible if $\mu|_F=0$. Iterating this procedure, we see that the restriction
to any subspace of $V$ (including $V$ itself) vanishes, hence $\mu=0$.  

\subsection{Schneider embedding}
\label{subsec_schneider_embedding}

The counterpart of Klain's embedding theorem for {\it odd} valuations was given
by Schneider. He showed in \cite{schneider96} that an odd, simple, continuous, translation invariant  valuation $\mu$ can be written as 
\begin{displaymath}
 \mu(K)=\int_{S(V)} f(v) dS_{n-1}(K,v),
\end{displaymath}
where $S_{n-1}(K,\cdot)$ is the $(n-1)$-th surface area measure of $K$ \cite{schneider_book93} and $f$ is an odd function on the unit sphere in $V$ (which will be denoted by $S(V)$). In particular, $\mu$ is of degree $n-1$.  

The function $f$ is  unique up to linear functions. Equivalently, $f$ is unique under the additional condition 
\begin{equation} \label{eq_orthogonality}
 \int_{S(V)} vf(v)dv=0. 
\end{equation}

Similarly as in the even case, this implies a description of odd valuations of a given degree. Namely, suppose $\mu \in \Val_k^-$. Then $\mu$ vanishes on $k$-dimensional sets, hence the restriction $\mu|_E$ to a $(k+1)$-dimensional subspace $E$ is simple and can be described by an odd function on the unit sphere of $E$ satisfying the condition \eqref{eq_orthogonality} with $V$ replaced by $E$. To put these functions into one object, one can use the {\it partial flag manifold} $\Flag_{k+1,1}$ consisting of pairs $(E,L)$, where $E \in \Gr_{k+1}(V)$ and $L$ is an oriented line in $E$. Then the {\it Schneider function} is an odd function on $\Flag_{k+1,1}$ (i.e. a function that changes sign if $(E,L)$ is replaced by $(E,-L)$). The space of continuous, odd functions on $\Flag_{k+1,1}$ is denoted by $C^{odd}(\Flag_{k+1,1} V)$. 

The valuation $\mu$ is uniquely determined by its Schneider function, as follows easily by induction on the dimension. Hence we get an embedding (the {\it Schneider embedding})
\begin{displaymath}
S : \Val_k^-  \hookrightarrow  C^{odd}(\Flag_{k+1,1} V).
\end{displaymath}

\subsection{Irreducibility theorem and smooth valuations}
\label{subsec_irr_thm}

Let $V$ be an $n$-dimensional vector space. Without fixing a Euclidean
structure on $V$, we still have the general linear group $\mathrm{GL}(V)$ acting on $V$
and on $\mathcal{K}(V)$. This action induces an action on $\Val$ by 
\begin{displaymath}
 g\mu(K):=\mu(g^{-1}K),
\end{displaymath}
which preserves degree and parity of a valuation. 

\begin{Theorem} {\bf (Alesker's irreducibility theorem)}\\ \label{thm_irr}
 The spaces $\Val_k^\epsilon, k=0,\ldots,n, \epsilon = \pm$ are irreducible
$\mathrm{GL}(V)$-representations. 
\end{Theorem}

We remind the reader that these spaces are in general infinite-dimensional
Banach spaces and that in this context, {\it irreducible} means that they do
not admit any non-trivial, invariant, closed subspaces. 

One way to understand the statement of the theorem is as follows. Start with a non-zero valuation $\mu \in \Val_k^\epsilon$ and consider its orbit under the group $\mathrm{GL}(V)$, i.e. the set of all $g\mu$. Then the space of linear combinations of such valuations are dense in $\Val_k^\epsilon$, which means that every valuation in $\Val_k^\epsilon$ may be approximated by these special ones.  

The proof of Theorem \ref{thm_irr} is contained in \cite{alesker_mcullenconj01}. It uses the Klain-Schneider embedding as well as heavy machinery from representation theory. 

Alesker's irreducibility theorem is of fundamental importance in algebraic
integral geometry. Let us explain the reason for this.

If we give some construction of translation invariant valuations, which does
not use any extra structure (like Euclidean metric), then we obtain a
$\mathrm{GL}(V)$-invariant subspace of $\Val$. By Alesker's irreducibility theorem, its
intersection with any of the spaces $\Val_k^\epsilon$ is either trivial or
dense. From this, one obtains several characterization theorems for
translation invariant valuations.  

We will see three main examples for this construction. The first is a positive
answer to a conjecture by McMullen \cite{mcmullen80}.  

\begin{Corollary} \label{cor_mcmullen}
 Valuations of the form $K \mapsto \vol(K+A)$, where $\vol$ is any Lebesgue
measure on $V$ and $A$ is a fixed convex body, span a dense subspace of
$\Val$. 
\end{Corollary}

The proof follows from the trivial observation that valuations of the form $K
\mapsto \vol(K+A)$ span a $\mathrm{GL}(V)$-invariant subspace of $\Val$ and that its
intersection with each $\Val_k^\epsilon$ is non-trivial.  

For the second example, we need the notion of {\it conormal cycle} of a compact
convex set. We suppose that $V$ is oriented and let $S^*V=V \times S(V^*)$ be the spherical cotangent bundle of
$V$, defined as follows. For $p \in V$, let $T_p^*V$ be the dual of the tangent space at $p$. On the space $T_p^*V \setminus \{0\}$, there is an equivalence relation given by $\xi \sim \xi'$ if and only if $\xi = \lambda \xi'$ for some real $\lambda>0$. The equivalence class of $\xi$ is denoted by $[\xi]$. 

The space $S^*V$ consists of the pairs $(p,[\xi])$, where $p \in V$, $\xi \in T^*_pV
\setminus \{0\}$. An element $(p,[\xi]) \in S^*V$ can be thought of as a pair
$(p,E)$, where $E=p+\ker \xi$ is an oriented affine hyperplane in $V$ containing $p$. 

The conormal cycle $N(K)$ of $K \in \mathcal{K}(V)$ is an
oriented $(n-1)$-dimensional Lipschitz submanifold in $S^*V$. It is given by the
set of all $(p,E)$ such that $p \in \partial K$ and $E$ is an oriented support
plane of $K$ at $p$.

If $V$ is a Euclidean vector space, then we can identify $S^*V$ with the sphere bundle $SV$. The image of the conormal cycle under this identification is the {\it normal cycle} of $K$. 

Let $\omega$ be a translation invariant $(n-1)$-form on $S^*V$ and $\phi$ be a
translation invariant $n$-form on $V$. Then the valuation 
\begin{equation} \label{eq_normal_cycle}
 K \mapsto \int_{N(K)} \omega+\int_K \phi
\end{equation}
is a continuous, translation invariant valuation. A valuation in $\Val$ of this
form is called {\it smooth} and $\Val^{sm}$ denotes the corresponding subspace.
By Alesker's irreducibility theorem, $\Val^{sm}$ is a dense subspace of $\Val$. 

The representation \eqref{eq_normal_cycle} opens the door to another central fact of algebraic integral geometry: {\it Smooth valuations can be extended to a large class of  compact non-convex sets}. Indeed, many compact sets $X \subset V$ admit a normal cycle $N(X)$ and \eqref{eq_normal_cycle} may be used to define $\mu(X)$. Examples of such sets are polyconvex sets (i.e. finite unions of convex sets), sets with positive reach, in particular smooth submanifolds (possibly with boundary or corner), compact sets which are definable in some o-minimal structure (see \cite{van_den_Dries98} for o-minimal structures and \cite{fu94, bernig_nc07} for the normal cycle of a definable set), in particular compact subanalytic or semialgebraic sets. 

Smooth valuations are natural from the viewpoint of
representation theory. Alesker's original definition of a smooth valuation in \cite{alesker03_un} was the following. 

\begin{Theorem} {\bf (Alesker, \cite{alesker_val_man1})}\\
A valuation $\mu \in
\Val$ is smooth if and only if the map
\begin{align*}
 \mathrm{GL}(V) & \to \Val\\
g & \mapsto g\mu 
\end{align*}
is smooth as a map from a Lie group to an infinite-dimensional Banach space.
\end{Theorem}

The proof uses the {\it Casselman-Wallach theorem} from representation theory. Furthermore, there is a natural way to endow $\Val^{sm}$ with a Fr\'echet space topology which is finer than the induced topology. 

The third application of Alesker's irreducibility theorem concerns {\it Crofton formulas} for even, homogeneous valuations. If $m$ is a (signed) translation invariant measure on the affine Grassmannian manifold $\AGr_{n-k}(V)$, then the valuation 
\begin{equation} \label{eq_crofton_measure}
 \mu(K):=\int_{\AGr_{n-k}(V)} \chi(K \cap E)dm(E)
\end{equation}
is an element of $\Val_k^+$. The signed measure $m$ is called {\it Crofton measure} of $\mu$. Since this construction is $\mathrm{GL}(V)$-invariant, it follows that the space of even, $k$-homogeneous valuations admitting such a Crofton measure is dense in $\Val_k^+$. If we restrict to smooth Crofton measures (i.e. measures which are given by integration over some smooth top-dimensional translation invariant form on $\AGr_{n-k}(V)$), then this subspace is precisely $\Val_k^{+,sm}$, i.e. the space of {\it smooth}, even, $k$-homogeneous valuations. For this last statement, the Casselman-Wallach theorem is used again, compare \cite{alesker_bernstein04}.

\subsection{Product}
\label{subsec_product}

One of the milestones of algebraic integral geometry is the introduction of a product structure on the space $\Val^{sm}$ by Alesker \cite{alesker04_product}. To define it, Alesker used his solution of McMullen's conjecture 
(Corollary \ref{cor_mcmullen}). If $A_1,A_2$ are smooth convex bodies with positive
curvature,
the valuations 
\begin{equation} \label{eq_def_phi_psi}
 \phi(K)=\vol_n(K+A_1), \psi(K)=\vol_n(K+A_2)
\end{equation}
are smooth and the Alesker product is defined by
\begin{equation} \label{eq_product}
 \phi \cdot \psi(K)=\vol_{2n}(\Delta K + A_1 \times A_2), 
\end{equation}
where $\Delta:V \to V \times V$ is the diagonal embedding. 

Alesker proved that this definition extends uniquely to a linear and continuous product 
\begin{align*}
 \Val^{sm} \times \Val^{sm} & \to \Val^{sm}\\
(\phi,\psi) & \mapsto \phi \cdot \psi. 
\end{align*}

If $\phi$ and $\psi$ are given as in \eqref{eq_def_phi_psi}, then 
\begin{align} \label{eq_phi_dot_psi}
 \phi \cdot \psi(K) & =\vol_{2n}(\Delta K + A_1 \times A_2) \nonumber \\
& = \int_V \int_V 1_{\Delta K+A_1 \times A_2}(x,y)dxdy \nonumber \\
& = \int_V \vol_n((y-A_2) \cap K + A_1)dy \nonumber \\
& = \int_V \phi((y-A_1) \cap K)dy.
\end{align}
This last expression extends to arbitrary $\phi \in \Val^{sm}$ by linearity. 

In the case of even smooth valuations, there is another description of the
product based on general Crofton formulas. We have seen that if $\phi \in
\Val_k^{+,sm}$, there is a smooth, translation invariant measure $m_\phi$ on the
space of $(n-k)$-planes in $V$ such that 
\begin{displaymath}
\phi(K)=\int_{\AGr_{n-k}(V)} \chi(K \cap E)dm_\phi(E).
\end{displaymath}

For $\psi$ as in \eqref{eq_def_phi_psi}, applying \eqref{eq_phi_dot_psi} and Fubini's theorem gives us 
\begin{align} \label{eq_prod_crofton}
 \phi \cdot \psi(K)& =\int_V \int_{\AGr_{n-k}(V)} \chi((y-A_2) \cap K \cap E)dm_\phi(E) dy \nonumber \\
& =\int_{\AGr_{n-k}(V)} \int_V \chi((y-A_2) \cap K \cap E) dy dm_\phi(E) \nonumber \\
& =\int_{\AGr_{n-k}(V)} \vol_n(K \cap E + A_2) dm_\phi(E) \nonumber \\
& =\int_{\AGr_{n-k}(V)} \psi(K \cap E) dm_\phi(E).
\end{align}
Again, this equation holds true for all $\psi \in \Val^{sm}$. In particular, we get that $\chi \cdot \psi=\psi$, i.e. {\it the Euler characteristic is the unit with respect to the Alesker product}.  

From \eqref{eq_prod_crofton} and \eqref{eq_int_vol_intersections} we see that the coefficient $c_k(K)$ in the general Hadwiger theorem \ref{thm_gen_hadwiger} is given by 
\begin{equation} \label{eq_coeff_gen_had}
 c_k(K)=\flag{n}{k}^{-1}\phi \cdot \mu_k(K), 
\end{equation}
which is already half of the ``algebraic'' proof of the general Hadwiger theorem. 

\subsection{Alesker-Poincar\'e duality} 

The Alesker product satisfies a remarkable Poincar\'e duality which is in fact a
central ingredient in the algebraic approach to kinematic formulas. By Klain's
theorem, the space $\Val_n$ (where $n$ is the dimension of $V$) is generated by
any Lebesgue measure. Fixing a Euclidean structure on $V$, we thus get an
isomorphism $\Val_n \cong \mathbb{C}$. Given two smooth valuations $\phi, \psi$,
let $\langle \phi,\psi\rangle \in \mathbb{C}$ be the image of the
$n$-homogeneous component of $\phi \cdot \psi$ under this isomorphism. 

Alesker proved that the pairing 
\begin{align*}
 \Val^{sm} \times \Val^{sm} & \to \mathbb{C}\\
 (\phi,\psi) & \mapsto \langle \phi,\psi\rangle
\end{align*}
is {\it perfect}, which means that the induced map 
\begin{displaymath}
 \mathrm{PD}:\Val^{sm} \to \Val^{sm,*}
\end{displaymath}
is injective and has dense image \cite{alesker04_product}. Roughly speaking, the space $\Val^{sm}$ is
{\it self-dual}. 

\subsection{Alesker-Fourier transform}

There is another remarkable duality on the space of translation invariant valuations, which shares many formal properties with the Fourier transform of functions. It was introduced by Alesker in the even case in \cite{alesker03_un} and in the odd case in \cite{alesker_fourier}. 

In invariant terms, the Alesker-Fourier transform is a map
\begin{displaymath}
 \wedge: \Val^{sm} \to \Val^{sm}(V^*) \otimes \Dens(V^*),
\end{displaymath}
where $\Dens(V^*)=\Lambda^n V \otimes \mathbb{C}$ denotes the $1$-dimensional space of complex-valued Lebesgue measures on $V^*$. Given a scalar product on $V$, we will identify $\Val^{sm}(V^*) \otimes \Dens(V^*)$ with $\Val^{sm}(V)$. 

The definition in the even case is easy to write down. If $\mu \in \Val_k^{sm,+}$ has Klain function $\Kl_\mu \in C^\infty(\Gr_k)$, then $\hat \mu \in \Val_{n-k}^{sm,+}$ is defined by the condition $\Kl_{\hat \mu}(E)=\Kl_\mu(E^\perp)$. Alesker showed that $\hat \mu$ indeed exists. One way to see this is to note that \eqref{eq_crofton_measure} can be rewritten in the form 
\begin{equation} \label{eq_crofton_measure2}
 \mu(K)=\int_{\Gr_kV} \vol(\pi_LK)dm(L),
\end{equation}
where $m$ is a (signed) smooth measure on $\Gr_kV$. Taking $m^\perp$ to be the image of $m$ under the map $L \mapsto L^\perp$, we can construct $\hat \mu$ by setting
\begin{displaymath}
 \hat \mu(K)=\int_{\Gr_{n-k}V} \vol(\pi_LK)dm^\perp(L).
\end{displaymath}
 
The definition in the odd case is much more involved and we refer to the original paper \cite{alesker_fourier} for the details. One of the main points of the construction is an odd version of a Crofton formula. If $\phi \in \Val_k^{-,sm}$, then one can write $\phi$ (non-uniquely) in the form
\begin{displaymath}
 \phi(K)=\int_{\Gr_{k+1}V} \psi_L(\pi_LK)dL,
\end{displaymath}
where $\psi_L \in \Val_k^{-,sm}(L)$ depends smoothly on $L$ and $\pi_L$ is the orthogonal projection onto $L$. The reader should compare this formula with \eqref{eq_crofton_measure2}. 

Using this formula, Alesker defined the Fourier transform on odd valuations in an inductive way and showed that the result does not depend on several choices (like the choice of the $\psi_L$).

The Alesker-Fourier transform satisfies a Plancherel-type formula: 
\begin{equation} \label{eq_plancherel}
 \hat{\hat \mu}(K)=\mu(-K).
\end{equation}

\subsection{Convolution}

Given a product and a Fourier transform, it is natural to consider the {\it
convolution} product on $\Val^{sm}$ in such a way that the Fourier transform is an algebra isomorphism between $(\Val^{sm},\cdot)$ and $(\Val^{sm},*)$, i.e.
\begin{equation} \label{eq_conv_prod}
 \widehat{\phi \cdot \psi}=\hat \phi * \hat \psi, \quad \phi,\psi \in \Val^{sm}.
\end{equation}

It was shown in \cite{bernig_fu06} that such a convolution product exists, and that there is the following equivalent definition, which is similar to \eqref{eq_product}. Suppose 
\begin{displaymath}
 \phi(K)=\vol(K+A_1), \psi(K)=\vol(K+A_2),
\end{displaymath}
where $A_1,A_2$ are smooth convex bodies with positive curvature. Then 
\begin{equation} \label{eq_convolution}
 (\phi * \psi)(K):=\vol(K+A_1+A_2),
\end{equation}
and the so-defined convolution extends to a unique linear and continuous product on 
$\Val^{sm}$. As the product, $*$ is commutative and associative. The volume is the unit in $(\Val^{sm},*)$. The degree of
$\phi * \psi$ is the sum of the degrees of $\phi$ and $\psi$ minus the
dimension $n$ of $V$. 

Note that the definition \eqref{eq_convolution} was given before the Alesker-Fourier transform was extended to odd valuations. Equation \eqref{eq_conv_prod}  in the odd case was established in \cite{alesker_fourier}. 

Let us make an important remark here. Since we use some volume in the
definition \eqref{eq_convolution} of the convolution, it is not independent of
the choice of a Euclidean scalar product on $V$. Without any choices, $*$ is
defined on the twisted space $\Val^{sm}(V^*) \otimes \Dens(V)$, where $\Dens(V) \cong \Lambda^n V^* \otimes \mathbb{C}$ denotes the $1$-dimensional space of complex-valued Lebesgue measures on $V$ . 

On the other hand, the product definition \eqref{eq_product} does not depend on any Euclidean structure, taking
$\vol_{2n}$ to be the product measure. The coordinate free version of the
Alesker-Fourier transform is an isomorphism
\begin{displaymath}
\Val^{sm}(V) \to \Val^{sm}(V^*) \otimes \Dens(V)
\end{displaymath}
and with these modifications \eqref{eq_conv_prod} is independent of a choice of Euclidean structure. 

Equation \eqref{eq_convolution} implies another nice property of the convolution. Namely, if $\phi,\psi$ are mixed volumes (see \cite{schneider_book93} for the definition and properties of mixed volumes)
 then their convolution product is again a mixed volume. More precisely, 
if $k +l \geq n$ and $A_1,\dots,A_{n-k}, B_1,\dots,B_{n-l}$ are convex bodies with smooth boundary and positive curvature, then the convolution product of the mixed volumes 
\begin{align*}
 \phi(K) & := V(K[k],A_1,\ldots,A_{n-k})\\ 
\psi(K) & := V(K[l],B_1,\ldots,B_{n-l}) 
\end{align*}
is the mixed volume 
\begin{displaymath}
\phi * \psi(K)= \binom{k+l}{k}^{-1} \binom{k + l}{n}V(K[k+l-n],A_1,\dots,A_{n-k},B_1,\dots,B_{n-l}).
\end{displaymath}

\subsection{Hard Lefschetz theorems}
\label{subsec_hlt}

The well-known {\it Hard Lefschetz theorem} from complex algebraic geometry states that iterates of the Lefschetz operator (multiplication by the symplectic form) realize the Poincar\'e-isomorphisms in the cohomology of K\"ahler manifolds. See  \cite{huybrechts05} for more information. 

In algebraic integral geometry, there is a similar theorem (in fact two versions of it). The Lefschetz operator is replaced by the multiplication with the first intrinsic volume $\mu_1$ (we fix some Euclidean structure here). The corresponding operator is denoted by 
\begin{displaymath}
 L:\Val^{sm}_* \to \Val^{sm}_{*+1}.
\end{displaymath}

Intertwining with the Alesker-Fourier transform, we get an operator 
\begin{displaymath}
 \Lambda:\Val^{sm}_* \to \Val^{sm}_{*-1}, \Lambda \phi=2\widehat{L \hat \phi}=2\mu_{n-1} * \phi.
\end{displaymath}
Explicitly, this operator is given by 
\begin{displaymath}
 \Lambda \mu(K)=\left.\frac{d}{dt}\right|_{t=0} \mu(K+tB),
\end{displaymath}
where $B$ is the unit ball and $K+tB$ is the parallel set of radius $t$ around $K$. 

From the Steiner formula \eqref{eq_steiner} and the trivial fact $\hat \mu_k=\mu_{n-k}$ one gets 
\begin{align*}
 L\mu_k & =\frac{(k+1)\omega_{k+1}}{2\omega_k} \mu_{k+1}\\
\Lambda \mu_k & = \frac{(n-k+1)\omega_{n-k+1}}{\omega_{n-k}} \mu_{k-1}.
\end{align*}

\begin{Theorem} \label{thm_hlt}
Let $V$ be an $n$-dimensional Euclidean vector space. 
\begin{enumerate}
 \item For $k \leq \frac{n}{2}$, the map 
\begin{displaymath}
 L^{n-2k}:\Val_k^{sm} \to \Val_{n-k}^{sm}
\end{displaymath}
is an isomorphism. 
\item For $k \geq \frac{n}{2}$, the map 
\begin{displaymath}
 \Lambda^{2k-n}:\Val_k^{sm} \to \Val_{n-k}^{sm}
\end{displaymath}
is an isomorphism. 
\end{enumerate}
\end{Theorem}

\begin{Corollary}
The multiplication operator 
 \begin{displaymath}
 L:\Val^{sm}_k \to \Val^{sm}_{k+1}
\end{displaymath}
is injective if $2k+1 \leq n$ and surjective if $2k+1 \geq n$. The derivation operator 
\begin{displaymath}
 \Lambda:\Val^{sm}_k \to \Val^{sm}_{k-1}
\end{displaymath}
is injective if $2k-1 \geq n$ and surjective if $2k-1 \leq n$.
\end{Corollary}

This corollary tells us that it is enough to understand valuations in the middle degree and that all other valuations are found by applying a simple operator to a valuation of middle degree. This is particularly useful when studying $G$-invariant valuations. The corollary also tells us that, roughly speaking, most valuations concentrate close to the middle degree. 

Several authors have contributed to the proof of Theorem \ref{thm_hlt}. Building on previous work with Bernstein \cite{alesker_bernstein04}, Alesker first proved both versions of the Hard Lefschetz theorem in the even case \cite{alesker03_un, alesker04}. The second version was extended to odd valuations in \cite{bernig_broecker07}. The proof used the Laplacian acting on differential forms on the sphere and some results from complex geometry (K\"ahler identities). Next, it was shown in \cite{bernig_fu06} that in the even case, both versions of the Hard Lefschetz theorem are in fact equivalent via the Alesker-Fourier transform (which was at that time defined only for even valuations). Finally, Alesker extended in \cite{alesker_fourier} the Fourier transform to odd valuations and derived the first version of the Hard Lefschetz theorem in the odd case from the second one.  

\section{Applications in integral geometry}
\label{sec_ftaig}

\subsection{Abstract Hadwiger-type theorem}

We have sketched in the first section how the kinematic formulas and Crofton
formulas can be easily proved with Hadwiger's theorem. A similar argument will give analogous
(although more complicated) formulas for all subgroups $G$ of $\mathrm{SO}(n)$ such that
$\dim \Val^G<\infty$. 

The next theorem was formulated by Alesker \cite{alesker_survey07}.
\begin{Theorem} \label{thm_abstract_hadwiger}
 A compact
subgroup $G$ of $\mathrm{SO}(n), n \geq 2$ satisfies $\dim \Val^G<\infty$ if
and only if $G$ acts transitively on the unit sphere. In this case, every $G$-invariant, translation invariant and continuous valuation is smooth. 
\end{Theorem}
 
Let us give the idea of the proof (taken from \cite{fu_survey06}). First of all,
remember
that smooth, translation invariant valuations are dense in the space of all
translation invariant valuations. A smooth valuation is given by integration
over the conormal cycle of some translation invariant differential form
$\omega$. If the valuation is $G$-invariant, then we may assume (by averaging
over the group) that $\omega$ is also $G$-invariant. If $G$ acts transitively
on the unit sphere, then a $G$-invariant, translation invariant differential
form on the sphere bundle is uniquely determined by its value at any given
point. Hence the space of all such forms is finite-dimensional. 

Now take any continuous, translation invariant, $G$-invariant valuation $\mu$. We may
approximate it by a sequence of smooth, translation invariant valuations.
Averaging these valuations with respect to the Haar measure on $G$, we may in
fact approximate $\mu$ by a sequence of smooth $G$-invariant, translation invariant
valuations. But this space is finite-dimensional, hence closed. Therefore $\mu$
itself belongs to this finite-dimensional space. In particular,  
\begin{displaymath}
 \Val^G \subset \Val^{sm}. 
\end{displaymath}

Let us now prove the inverse implication. The Klain embedding \eqref{eq_klain_embedding} in the case $k=1$ induces an isomorphism 
\begin{displaymath}
 \Kl:\Val_1^+ \cong C^\infty(\Gr_1V).
\end{displaymath}
This follows from the fact that the {\it cosine transform} is an isomorphism on even smooth functions on the unit sphere \cite{klain_rota}. Since $\mu$ is $G$-invariant if and only if $\Kl_\mu$ is $G$-invariant, we have 
\begin{displaymath}
 \Kl:\Val_1^{G,+} \cong C^\infty(\Gr_1V)^G.
\end{displaymath}
If $G$ does not act transitively on the sphere, then the space of smooth  $G$-invariant functions on the projective space $\Gr_1V=\mathbb{P}V$ is infinite-dimensional. Therefore $\Val_1^{G,+}$ is also infinite-dimensional, which implies that $\Val^G$ is infinite-dimensional. 

The classification of connected compact Lie groups $G$ acting transitively on some
sphere is a topological problem which was solved by Montgomery-Samelson \cite{montgomery_samelson43} and Borel \cite{borel49}. There are six infinite lists
\begin{equation} \label{eq_transitive_groups}
 \mathrm{SO}(n), \mathrm{U}(n), \mathrm{SU}(n), \mathrm{Sp}(n), \mathrm{Sp}(n) \cdot \mathrm{U}(1), \mathrm{Sp}(n)\cdot \mathrm{Sp}(1)
\end{equation}
and three exceptional groups
\begin{equation} \label{eq_exceptional_groups}
 \mathrm{G}_2, \mathrm{Spin}(7), \mathrm{Spin}(9). 
\end{equation}

These groups are important in differential geometry and topology, since the holonomy group of an irreducible non-symmetric Riemannian manifold is always from this list and each group from this list except $\mathrm{Sp}(n) \cdot \mathrm{U}(1)$ and $\mathrm{Spin}(9)$ does appear as the holonomy group of such a manifold.

There are various natural inclusions among these groups:  
\begin{multline*}
 \mathrm{U}(n), \mathrm{SU}(n) \subset \mathrm{SO}(2n), \quad \mathrm{Sp}(n), \mathrm{Sp}(n)\cdot \mathrm{U}(1),
\mathrm{Sp}(n)\cdot \mathrm{Sp}(1) \subset \mathrm{SO}(4n), \\ 
\mathrm{SU}(4) \subset \mathrm{Spin}(7), \mathrm{G}_2 \subset  \mathrm{SO}(7), \mathrm{SU}(3) \subset \mathrm{G}_2 \subset  \mathrm{Spin}(7) \subset \mathrm{SO}(8), \mathrm{Spin}(9) \subset \mathrm{SO}(16).
\end{multline*}
The last two inclusions are the spin representations. We refer to \cite{besse_book} for more information on holonomy groups.

\subsection{The kinematic coproduct}

The first thing we need to do in order to relate the kinematic formulas to the
algebraic structures introduced in the previous section is to give a more
abstract description of these formulas. 

Let $V$ be a Euclidean vector space and let $G$ be a subgroup of $\mathrm{SO}(V)$ which
acts transitively on the unit sphere. We have seen that in this case, the space $\Val^G$ is finite-dimensional and consists only of smooth valuations. 

If $\phi_1,\ldots,\phi_m$ is a basis of $\Val^G$, then by the same trick as in
Subsection \ref{subsec_hadwiger} we obtain kinematic formulas 
\begin{equation} \label{eq_kf_g}
 \int_{\bar G} \phi_i(K \cap \bar g L)d\bar g=\sum_{k,l=1}^m c_{k,l}^i
\phi_k(K)\phi_l(L). 
\end{equation}

There is a very nice and clever way to encode these formulas in a purely
algebraic way. For this, Fu \cite{fu06} defined the kinematic operator 
\begin{align*}
 k_G:\Val^G & \to \Val^G \otimes \Val^G\\
\phi_i & \mapsto \sum_{k,l=1}^m c_{k,l}^i \phi_k \otimes \phi_l.
\end{align*}

This map is in fact a {\it cocommutative, coassociative coproduct} on $\Val^G$.
Let us remind the reader of the definition of a coproduct. Loosely speaking, we
write down the corresponding usual property (commutativity or associativity) in
terms of a commuting diagram and reverse all arrows to obtain the co-property. 

For instance, cocommutativity means that the
following diagram commutes:

\begin{displaymath}
\xymatrix{\Val^G \ar[r]^<<<<<{k_G} \ar[d]_{id} & \Val^G \otimes \Val^G \ar[d]_{\iota}
\\ \Val^G \ar[r]^<<<<<{k_G} & \Val^G \otimes \Val^G.}
\end{displaymath}
Here $\iota$ is the map that interchanges the factors of $\Val^G \otimes
\Val^G$. 

In more concrete terms, this says that the
coefficients in the kinematic formula \eqref{eq_kf_g} satisfy $c_{k,l}^i=c_{l,k}^i$, which
expresses the symmetry of the formula (in $K$ and $L$) as in Subsection
\ref{subsec_kinformula_son}. 

The coassociativity property is the commutativity of the following diagram: 
\begin{displaymath}
 \xymatrix{\Val^G \ar[r]^<<<<<{(k_G
\otimes id) \circ k_G} \ar[d]_{id} & \Val^G \otimes \Val^G \otimes \Val^G 
\ar[d]_{id \otimes id \otimes id}\\
\Val^G \ar[r]^<<<<<{(id \otimes k_G) \circ k_G} &
\Val^G \otimes \Val^G \otimes \Val^G.}
\end{displaymath}

This property is equivalent to the formula 
\begin{displaymath}
\sum_r c_{r,m}^ic_{k,l}^r=\sum_r c_{r,l}^ic_{k,m}^r, 
\end{displaymath}
and this comes just from Fubini's theorem, compare Subsection
\ref{subsec_kinformula_son}. 
 
In a similar vein, there are additive kinematic formulas for $G$:
\begin{equation} \label{eq_add_kf_g}
 \int_{G} \phi_i(K + g L)d\bar g=\sum_{k,l=1}^m d_{k,l}^i
\phi_k(K)\phi_l(L). 
\end{equation}

which can be encoded by the cocommutative, coassociative coproduct 
\begin{align} \label{eq_a_G}
a_G:\Val^G & \to \Val^G \otimes \Val^G \nonumber \\
\phi_i & \mapsto \sum_{k,l=1}^m d_{k,l}^i \phi_k \otimes \phi_l.
\end{align}

\subsection{Fundamental theorem of algebraic integral geometry}
\label{subsec_ftaig}

The {\it fundamental theorem of algebraic integral geometry} relates the
kinematic coproduct and the product structure and is the basis for a fuller
understanding of the kinematic formulas \eqref{eq_kf_g}. 

\begin{Theorem} \label{thm_ftaig}
Let $G$ be a group acting transitively on the unit sphere, $m_G:\Val^G \otimes
\Val^G \to \Val^G$ the restriction of the Alesker product to $\Val^G$;
$\mathrm{PD}_G:\Val^G \to \Val^{G*}$ the restriction of the Alesker-Poincar\'e duality to
$\Val^G$ and $k_G$ the kinematic coproduct. Then the following diagram commutes
\begin{displaymath}
 \xymatrix{\Val^G \ar[r]^<<<<<<{k_G} \ar[d]_{\mathrm{PD}_G} & \Val^G \otimes \Val^G
\ar[d]_{\mathrm{PD}_G \otimes \mathrm{PD}_G} \\
\Val^{G*} \ar[r]^<<<<<{m_G^*} & \Val^{G*} \otimes \Val^{G*}.}
\end{displaymath}
\end{Theorem}

This theorem, based on a basic version which we discuss below, was proven in \cite{bernig_fu06}. 

Let us work out the most important case, namely the principal kinematic formula $k_G(\chi)$. First note that $\Val^G \otimes \Val^G=\Hom(\Val^{G*},\Val^G)$. We may thus regard $k_G(\chi)$ as a map from $\Val^{G*}$ to $\Val^G$. Recall that $\mathrm{PD}_G$ is a map from $\Val^G$ to $\Val^{G*}$. 

Given $K \in \mathcal{K}(V)$, let $\tau_K \in \Val^{G*}$ be defined by $\tau_K(\phi):=\phi(K)$. Then the $\tau_K$ span $\Val^{G*}$ and we get  
\begin{displaymath}
 k_G(\chi)(\tau_K)(\cdot)=\int_{\bar G} \chi(K \cap \bar g \cdot) d\bar g \in \Val^G
\end{displaymath}
and hence for any $\phi \in \Val^G$ by \eqref{eq_prod_crofton} 
\begin{displaymath}
 \phi \cdot k_G(\chi)(\tau_K)(\cdot)=\int_{\bar G} \phi(K \cap \bar g \cdot) d\bar g \in \Val^G.
\end{displaymath}
We plug in a ball $B_R$ of radius $R$ into this equation. If $R$ is large, the measure of all $\bar g$ with $K \subset \bar g B_R$ is approximately $\vol(B_R)$, while the measure of all $\bar g$ with $K \cap \bar g B_R \neq \emptyset,K$ is $o(R^n)$. It follows that the $n$-th homogeneous component of $\phi \cdot k_G(\chi)(\tau_K)$ is given by $\phi(K)\vol$. In other words, 
\begin{displaymath}
 \mathrm{PD}_G(k_G(\chi)(\tau_K))(\phi)=\phi(K)=\tau_K(\phi),
\end{displaymath}
which implies that 
\begin{displaymath}
 \mathrm{PD}_G \circ k_G(\chi)=Id. 
\end{displaymath}
Hence {\it $\mathrm{PD}_G$ and $k_G(\chi)$ are inverse to each other}, and this statement is equivalent to the fact that 
\begin{displaymath}
 (\mathrm{PD}_G \otimes \mathrm{PD}_G) \circ k_G(\chi)=m_G^* \circ \mathrm{PD}_G(\chi),
\end{displaymath}
which follows from Theorem \ref{thm_ftaig}. See also  \cite{fu06, bernig_fu06} for more details. 
 
The fundamental theorem of algebraic integral geometry says roughly that the knowledge of $k_G$ is the same as
the knowledge of $m_G$. It may be used in two ways. If we know $k_G$, then we
may first compute $\mathrm{PD}_G:=k_G(\chi)^{-1}$ and, using the above diagram, we may
compute the whole product structure. Conversely, knowing the product, we can
compute $\mathrm{PD}_G$ and hence $k_G$. This is how the theorem will be used in the
sequel. 

Nevertheless, in concrete situations, things turn out to be not so easy, since
in order to compute $m_G^*$ we have to invert some potentially huge matrix which
might be a challenge. We will come back to this point when we discuss the
hermitian case. 

A consequence from the fundamental theorem of algebraic integral geometry is
\begin{equation}
 k_G(\phi \cdot \psi)=(\phi \otimes \chi) \cdot k_G(\psi)=(\chi \otimes \psi) \cdot k_G(\phi), \quad \phi,\psi \in \Val^G.   
\end{equation}

We give a proof of the more general statement
\begin{equation} \label{eq_prod_coprod}
 \int_{\bar G} \phi \cdot \psi(K \cap \bar g L)d\bar g =((\phi \otimes \chi) \cdot k_G(\psi))(K,L),
\end{equation}
where $\psi$ is supposed to be smooth and translation invariant, but not necessarily $G$-invariant. 

By linearity and density, it is enough to assume that $\phi$ has the form $\phi(K)=\vol(K+A)$ for some smooth convex body $A$ with positive curvature.
Then 
\begin{displaymath}
 \phi \cdot \psi(K \cap \bar g L)=\int_V \psi((x-A) \cap K \cap \bar g L)dx
\end{displaymath}
by \eqref{eq_phi_dot_psi} and hence 

\begin{align*}
 \int_{\bar G} \phi \cdot \psi(K \cap \bar g L)d\bar g & = \int_{\bar G} \int_V \psi((x-A) \cap K \cap \bar g L)dx d\bar g \\
& = \int_V \int_{\bar G}  \psi((x-A) \cap K \cap \bar g L) d\bar g dx\\ 
& = \int_V k_G(\psi)((x-A) \cap K,L) dx\\ 
& = (\phi \otimes \chi) \cdot k_G(\psi)(K,L).
\end{align*}

In the special case $G=\mathrm{SO}(n)$, $\psi=\chi$, Equation \eqref{eq_prod_coprod} and the principal kinematic formula \eqref{eq_kin_formula_son} imply the general Hadwiger theorem \ref{thm_gen_hadwiger}: 
\begin{align*}
 \int_{\overline{\mathrm{SO}(n)}} \phi (K \cap \bar g L)d\bar g & =((\phi \otimes \chi) \cdot k_{\mathrm{SO}(n)}(\chi))(K,L)\\
& = \sum_{k=0}^n \flag{n}{k}^{-1} (\phi \cdot \mu_k)(K) \mu_{n-k}(L)\\
& = \sum_{k=0}^n c_k(K) \mu_{n-k}(L),
\end{align*}
where 
\begin{displaymath}
 c_k(K)=\flag{n}{k}^{-1} (\phi \cdot \mu_k)(K)=\int_{\AGr_{n-k}(V)} \phi(K \cap E)dE
\end{displaymath}
by \eqref{eq_coeff_gen_had}. 

In the same situation, Nijenhuis' observation becomes evident. We set $t:=\frac{2}{\pi}\mu_1 \in \Val_1^{\mathrm{SO}(n)}$. By the Hard Lefschetz theorem \ref{thm_hlt} we must have $t^n=c \vol_n$ for some constant $c \neq 0$. Therefore   
\begin{displaymath}
 \Val^{\mathrm{SO}(n)}=\mathbb{C}[t]/(t^{n+1}).
\end{displaymath}

Then we have $\mathrm{PD}(t^i)=c (t^{n-i})^*$, where $\left\{(t^k)^*, k=0,\ldots,n\right\}$ is the dual basis to the basis $\left\{t^k, k=0,\ldots,n\right\}$ of $\Val^{\mathrm{SO}(n)}$. From Theorem \ref{thm_ftaig} it follows that  
\begin{displaymath}
 k_G(t^i)=\frac{1}{c} \sum_{k+l=n+i} t^k \otimes t^l.
\end{displaymath}

Setting $\tilde \mu_k=\frac{1}{c}t^k$ thus gives us 
\begin{displaymath}
 k_G(\tilde \mu_i)=\sum_{k+l=n+i} \tilde \mu_k \otimes \tilde \mu_l.
\end{displaymath}

In fact, it is easily computed (see Subsection \ref{subsec_hlt} or \cite{bernig_fu_hig}) that
\begin{displaymath}
 t^k=\frac{k!\omega_k}{\pi^k}\mu_k,
\end{displaymath}
hence $c=\frac{n!\omega_n}{\pi^n}$ and 
\begin{displaymath}
 \tilde \mu_k=\frac{\pi^n k! \omega^k}{\pi^k n! \omega_n} \mu_k.
\end{displaymath}

\subsection{Additive formulas}

There is a similar statement relating the convolution product to the additive
kinematic formulas \eqref{eq_add_kf_g}. It was proved in \cite{bernig_fu06} under the assumption $\Val^G \subset \Val^+$, which turns out to be always the case \cite{bernig_g2}. 
 
\begin{Theorem}
Let $G$ be a group acting transitively on the unit sphere. Let $a_G$ be the additive kinematic coproduct, see \eqref{eq_a_G}. Let
$c_G:\Val^G \otimes \Val^G \to \Val^G$ be the restriction of the convolution to
$\Val^G$. Then the following diagram commutes
\begin{displaymath}
 \xymatrix{\Val^G \ar[r]^<<<<<{a_G} \ar[d]_{\mathrm{PD}_G} & \Val^G \otimes \Val^G
\ar[d]_{\mathrm{PD}_G \otimes \mathrm{PD}_G} \\
\Val^{G*} \ar[r]^<<<<<{c_G^*} & \Val^{G*} \otimes \Val^{G*}.}
\end{displaymath}
\end{Theorem}

\begin{Corollary}
Kinematic formulas \eqref{eq_kf_g} and additive kinematic formulas \eqref{eq_add_kf} are related by the formula 
\begin{equation}
 a_G=(\wedge \otimes \wedge) \circ k_G \circ \wedge.
\end{equation}
Explicitly, this means that if the kinematic formulas are given by 
\begin{displaymath}
 \int_{\bar G} \phi_i(K \cap \bar g L)d\bar g=\sum_{k,l=1}^m c_{k,l}^i
\phi_k(K)\phi_l(L),
\end{displaymath}
in some basis $\phi_1,\ldots,\phi_m$ of $\Val^G$, then the additive kinematic formulas in the Fourier-dual basis $\hat \phi_1,\ldots,\hat \phi_m$ are given by 
\begin{displaymath}
 \int_G \hat \phi_i(K+g L)dg=\sum_{k,l=1}^m c_{k,l}^i \hat \phi_k(K)\hat
\phi_l(L),
\end{displaymath}
with the same constants. 
\end{Corollary}
 
This corollary gives a precise meaning to the fact which we have mentioned in
Subsection \ref{subsec_kinformula_son}: {\it Kinematic formulas and additive kinematic formulas are
dual to each other.} 

This explains also the observation from Subsection \ref{subsec_kinformula_son}: since in some basis of $\Val^{\mathrm{SO}(n)}$ all coefficients of $k_{\mathrm{SO}(n)}$ are $1$, the same holds true for $a_{\mathrm{SO}(n)}$ in the Fourier-dual basis.

\section{The hermitian case}
\label{sec_un}

In his 1976 book on integral geometry \cite{santalo76}, Santal\'o wrote that {\it Integral geometry on complex spaces has not been sufficiently developed and probably deserves further study.}

In the previous two sections, we have described the theoretical
framework relating algebraic structures on valuations and
integral-geometric formulas. The aim of this section is to show how this works
in practice for the first non-classical example from list
\eqref{eq_transitive_groups}, namely the group $G=\mathrm{U}(n)$.  

We let $V \cong \mathbb{C}^n$ be a complex vector space of (complex)
dimension $n$, endowed with a hermitian inner product $H$. Recall that $H$ is 
\begin{enumerate}
 \item conjugate linear in the first component and linear in the second component, i.e.
\begin{displaymath}
 H(\lambda v,\mu w)=\bar \lambda H(v,w) \mu, \quad v,w \in V, \lambda, \mu \in \mathbb{C},
\end{displaymath}
\item conjugate symmetric, i.e. $H(w,v)=\overline{H(v,w)}$ and
\item positive definite, i.e. $H(v,v)>0$ for $v \neq 0$. 
\end{enumerate}

The subgroup of $\mathrm{GL}(V,\mathbb{C})$ fixing $H$ is the unitary group
$\mathrm{U}(n)$. 

The real part of $H$ is a real inner product on $V$, while the imaginary part of $H$ is a symplectic form $\Omega$ on $V$. In particular, $\mathrm{U}(n)$ is a subgroup of $\mathrm{SO}(2n)$. 

Before going into details, let us remark that $-1 \in \mathrm{U}(n)$, hence all
unitarily invariant valuations are even. 

\subsection{$\Val^{\mathrm{U}(n)}$ as a vector space}
\label{subsec_valun_vec}

The abstract Hadwiger theorem \ref{thm_abstract_hadwiger} tells us that $\dim
\Val^{\mathrm{U}(n)}<\infty$, but it says nothing about the actual value of this
dimension. Alesker showed in \cite{alesker_mcullenconj01} that 
\begin{equation} \label{eq_dim_un}
 \dim \Val_k^{\mathrm{U}(n)}=\min\left\{\left\lfloor
\frac{k}{2}\right\rfloor,\left\lfloor\frac{2n-k}{2}\right\rfloor\right\}+1. 
\end{equation}

Note that these dimensions have the typical behavior predicted by the Hard
Lefschetz Theorem \ref{thm_hlt}: they are increasing for degrees smaller than
half the (real) dimension and decreasing for degrees bigger than half the
(real) dimension.  

There are various ways of proving this formula. Alesker's original proof used
representation theoretical methods to decompose the space of even valuations on
an $2n$-dimensional vector space as a direct sum of irreducible
$\mathrm{SO}(2n)$-modules. Since it is known which irreducible $\mathrm{SO}(2n)$-modules contain
a $\mathrm{U}(n)$-invariant vector, the above formula follows easily. 

A second possible proof goes as follows. Since $\Val^{\mathrm{U}(n)} \subset \Val^{sm}$,
we can represent each unitarily invariant valuation by a pair $(\omega,\phi)$ of
differential forms as in \eqref{eq_normal_cycle}. Since we may average over the
group, we may actually take $\omega, \phi$ to be $\mathrm{U}(n)$-invariant too. But the
$\mathrm{U}(n)$-invariant, translation invariant smooth forms on the sphere bundle $SV$
can be explicitly described. This was carried out by Park \cite{park02} using
the
{\it first fundamental theorem} for the group $\mathrm{U}(n)$. Different pairs
$(\omega,\phi)$ may induce the same valuation. Fortunately, one can
characterize the kernel of the normal cycle map in terms of a certain
second-order differential operator, called {\it Rumin operator} which
was introduced by Rumin in \cite{rumin94}. This works even in the much
more general setting of {\it valuations on manifolds}, see
\cite{bernig_broecker07}. The Rumin operator of the unitarily and
translation invariant forms on $SV$ was (somehow implicitly) computed in
\cite{bernig_fu_hig}. These computations imply \eqref{eq_dim_un}.  

A third way to prove \eqref{eq_dim_un} is sketched in \cite{bernig_qig}. It uses the fact that $\Val_k^G$ and some spaces of $G$-invariant differential forms on the unit sphere bundle $SV$ fit into an exact sequence.

Knowing the dimension of $\Val^{\mathrm{U}(n)}$, the next question is to find a basis.
Alesker gave in fact two of them, which are dual to each other with respect to
the Alesker-Fourier transform. The idea is to mimic the definition of the
intrinsic volumes in \eqref{eq_int_vol_projections} and
\eqref{eq_int_vol_intersections} and using complex Grassmannians instead
of real ones. Using intersections with complex planes, Alesker defined 
\begin{equation*}
U_{k,p}(K) := \int_{\overline{\Gr}^\mathbb{C}_{n-p}} \mu_{k-2p}(  K\cap \bar E) \,
d\bar E.
\end{equation*}

The $U_{k,p}$, as $p$ ranges over $0,1,\ldots,\min \left\{\left\lfloor \frac k
2\right\rfloor, \left\lfloor \frac {2n-k}2 \right\rfloor\right\}$, constitute a
basis of $\Val_k^{\mathrm{U}(n)}$. 

Fu renormalized these valuations by setting 
\begin{align*}
 t & := \frac{2}{\pi} \mu_1=\frac{2}{\pi} U_{1,0} \in \Val_1^{\mathrm{U}(n)}\\
s & := nU_{2,1} \in \Val_2^{\mathrm{U}(n)}
\end{align*}
which implies that 
\begin{displaymath}
s^pt^{k-2p} = \frac{(k-2p)!n!\omega_{k-2p}}{(n-p)!\pi^ {k-2p}} U_{k,p}. 
\end{displaymath}

The second basis given by Alesker uses projections onto complex
subspaces instead of intersections:  
\begin{equation*}
C_{k,q}(K) := \int_{\Gr^\mathbb{C}_q} \mu_k(\pi_E(K)) \, dE.
\end{equation*}

As $q$ ranges over all values from $n-\min \left\{\left\lfloor \frac k
2\right\rfloor, \left\lfloor \frac {2n-k}2 \right\rfloor\right\}$ to $n$, the $C_{k,q}$ constitute a basis of
$\Val_k^{\mathrm{U}(n)}$. Up to a normalizing constant, the Fourier transform of
$U_{k,p}$ is $C_{2n-k,n-p}$. 

\subsection{$\Val^{\mathrm{U}(n)}$ as an algebra}

The monomials $s^pt^{k-2p}$, with $0 \leq k \leq 2n$ and $0 \leq p \leq \min
\left\{\left\lfloor \frac k2\right\rfloor, \left\lfloor \frac {2n-k}2
\right\rfloor\right\}$, constitute a
basis of $\Val^{\mathrm{U}(n)}$. We therefore speak of the {\it monomial basis} or
the {\it $ts$-basis} of $\Val^{\mathrm{U}(n)}$.  

We have a graded algebra epimorphism
\begin{displaymath}
 \mathbb{C}[t,s] \twoheadrightarrow \Val^{\mathrm{U}(n)},
\end{displaymath}
where $t,s$ on the left hand side are formal variables of degree $1$ resp.
$2$ (in the following, the distinction between variables and actual valuations
will not be made, which is quite in the spirit of algebraic integral geometry).
The kernel of this map is an ideal $I_n$ in $\mathbb{C}[t,s]$, which, by Hilbert
basis theorem, must be generated by finitely many polynomials.  

There is a relatively easy way to compute these polynomials, which was given by Fu \cite{fu_les_diablerets}. 

First,
one deduces from \eqref{eq_dim_un} that $I_n$ is generated by two polynomials
$f_{n+1}$ and $f_{n+2}$ of total degree $n+1$ and $n+2$ respectively.  

Next, by Alesker-Poincar\'e duality, in order to show that some polynomial $f$
of total degree $d$ in $t$ and $s$ is zero, it is enough to show that $f \cdot
s^pt^{2n-d-2p}=0$ for all $p$. Since $\Val_{2n}^{\mathrm{U}(n)}$ is spanned by the
Lebesgue measure, this amounts to some combinatorial identity among the
coefficients of $f$ once we know how to evaluate the monomials $s^pt^{2n-2p}$ on
a unit ball. Using the {\it transfer principle}, which relates valuations on $\mathbb{C}^n$ and on $\mathbb{CP}^n$, \cite{howard93}, one can compute
these values. The final result (which was first
proved by Fu in \cite{fu06} using another method) is as follows: 

\begin{Theorem} \label{thm_fu}
 There is an isomorphism between graded algebras 
\begin{displaymath}
 \Val^{\mathrm{U}(n)} \cong \mathbb{C}[t,s]/(f_{n+1},f_{n+2}),
\end{displaymath}
where 
\begin{displaymath}
 \log(1+t+s)=f_1+f_2+f_3+\ldots=t+\left(s-\frac{t^2}{2}\right)+\left(-st+\frac{
t^3}{3}\right)+\ldots
\end{displaymath}
is the expansion in (weighted) homogeneous polynomials. 
\end{Theorem}

As explained in Subsection \ref{subsec_ftaig}, from the product structure, we
can compute $\mathrm{PD}_{\mathrm{U}(n)}$ and $m_{\mathrm{U}(n)}$ and therefore $k_{\mathrm{U}(n)}$. Theorem
\ref{thm_fu} thus implies the knowledge of the kinematic formulas in the
$ts$-basis. 

Working this out in higher dimensions is rather cumbersome, because some huge
matrix has to be inverted. Also, one would like to have not only the
value of the coefficients in the kinematic formulas \eqref{eq_kf_g}, but some closed forms. They seem to be hard to obtain from Theorem \ref{thm_fu}. Another missing point is
the knowledge of the kinematic formula in another basis of $\Val^{\mathrm{U}(n)}$, for
instance in the $C$-basis.  
 
\subsection{Hermitian intrinsic volumes and Tasaki valuations}

It seems difficult to describe the value of a basis element of the $ts$-basis
on, say a polytope or a submanifold (since all unitarily invariant valuations are smooth, they may be canonically extended to submanifolds with boundary or corners, see Subsection \ref{subsec_irr_thm}). Therefore we introduce another, more geometric basis. This mimics the third
characterization of the intrinsic volumes in Subsection \ref{subsec_int_vol}.  

Recall that a real subspace $E$ of $V$ is called {\it isotropic} if the
restriction of the symplectic form to $E$ vanishes. Then the dimension of $E$
does not exceed $n$, and an isotropic subspace of dimension $n$ is called {\it
Lagrangian}. We call $E$ of type $(k,q)$ if $E$ can be written as the
orthogonal sum of a complex subspace of (complex) dimension $q$ and an isotropic
subspace of dimension $k-2q$. Then $k-q \leq n$. 

\begin{Theorem}
 There is a unique valuation $\mu_{k,q} \in \Val_k^{\mathrm{U}(n)}$ whose Klain function
evaluated at a subspace of type $(k,q')$ equals $\delta_{qq'}$. Moreover, 
\begin{displaymath}
 \hat \mu_{k,q}=\mu_{2n-k,n-k+q}.
\end{displaymath}
\end{Theorem}

The idea of the construction of $\mu_{k,q}$ is as follows. We know from
the discussion in Subsection \ref{subsec_valun_vec} that every unitarily
invariant valuation of degree $k<2n$ is given by integration over the normal
cycle of some translation invariant, unitarily invariant differential form on
$SV$. Park \cite{park02} showed that the algebra of these forms is generated by
three $1$-forms and four $2$-forms, and integrating a suitable product of these
basic forms over the normal cycle yields the valuation $\mu_{k,q}$. 

Since the $\mu_{k,q}$ with $\max(0,k-n)\leq q \leq \lfloor \frac k
2\rfloor$ are linearily independent, it follows from
\eqref{eq_dim_un} that they form a basis of $\Val_k^{\mathrm{U}(n)}$. 

Finally, the statement on the Fourier transform boils down to the fact that the
orthogonal complement of a subspace of type $(k,q)$ is of type $(2n-k,n-k+q)$, which
is easy to prove.

A version of these valuations was considered by Tasaki. He showed
that the orbits of the $\mathrm{U}(n)$-action on $\Gr_k(V)$ are characterized by
$\left\lfloor \frac{\min\{k,2n-k\}}{2}\right\rfloor$ {\it K\"ahler angles}. We use a slight modification of his construction. Let $p:=\lfloor \frac{k}{2}\rfloor$.  Given a $k$-dimensional subspace $E \subset V$, the restriction of the symplectic form $\Omega$ of $V$ to $E$ can be written as 
\begin{displaymath}
 \Omega|_E=\sum_{i=1}^p \cos \theta_i \alpha_{2i-1} \wedge \alpha_{2i},
\end{displaymath}
where $\alpha_1,\ldots,\alpha_k$ is dual to an orthonormal basis of $E$ and $0 \leq \theta_1 \leq \ldots \leq \theta_p \leq \frac{\pi}{2}$. The $p$-tuple $(\theta_1,\ldots, \theta_p)$ is called {\it multiple K\"ahler angle} of $E$. 

For instance, a subspace is isotropic if all its K\"ahler angles are
$\frac{\pi}{2}$, while it is complex if all K\"ahler angles are $0$. More
generally, a subspace is of type $(k,q)$ if $q$ of its K\"ahler angles are $0$ and
the remaining $p-q$ K\"ahler angles are $\frac{\pi}{2}$. Tasaki \cite{tasaki00} showed that two $k$-dimensional subspaces belong to the same $\mathrm{U}(n)$-orbit if and only if their multiple K\"ahler angles agree. 

The {\it Tasaki valuations} $\tau_{k,q} \in \Val^{\mathrm{U}(n)}, 0 \leq q \leq
p$ are defined by their Klain function: 
\begin{equation}
\Kl_{\tau_{k,q}}(E) =\sigma_q(\cos^2\theta_1(E),\dots,\cos^2\theta_p(E)),
\end{equation}
where $\sigma_q$ is the the $q$th elementary symmetric function. 

It is of course elementary to compute the relations between the Tasaki
valuations and the hermitian intrinsic volumes: 
\begin{equation} \label{eq_relation_tasaki_hiv}
\tau_{k,q}=\sum_{i=q}^{\lfloor k/2\rfloor} \binom{i}{q} \mu_{k,i}, \quad
\mu_{k,q}=\sum_{i=q}^{\lfloor k/2\rfloor} (-1)^{i-q} \binom{i}{q} \tau_{k,i}.
\end{equation}

If $M$ is a compact $k$-dimensional manifold, then the canonical extension of $\tau_{k,q}$ to $M$ is given by 
\begin{displaymath}
 \int_M \sigma_q (\cos^2 \Theta(T_x M)) \, dx.
\end{displaymath}
Using such expressions, Tasaki \cite{tasaki03} formulated general {\it Poincar\'e formulas}, which are special instances of the principal kinematic formula $k_{\mathrm{U}(n)}(\chi)$, with $K,L$ replaced by compact submanifolds of complementary dimension. 

\subsection{Kinematic formulas}
\label{subsec_kun}

Let us now explain, in an informal style, how the hermitian intrinsic
volumes may be used to compute the relations between the different bases
($U$-basis and $C$-basis), and to compute the kinematic formulas. 

One can easily compute the derivation operator $\Lambda$ (compare
Subsection \ref{subsec_hlt}) on the hermitian intrinsic volumes. This comes from the fact that
the hermitian intrinsic volumes are given by integration over the normal cycle
of certain differential forms. The operator $\Lambda$ corresponds to a certain Lie
derivative on the level of forms which is easy to compute. 

Since we also know the Alesker-Fourier transform of $\mu_{k,q}$, we can compute
$L\mu_{k,q}$ (which is multiplication by $t$, up to a factor). Now a crucial
(and somehow mysterious) observation is that (some renormalizations of) $L$ and
$\Lambda$ and some degree counting operator define a
representation of the Lie algebra $\mathfrak{sl}_2$ on $\Val^{\mathrm{U}(n)}$. In the general translation invariant setting $\Val$, this is not the case. 

The next observation is that $\mu_{n,0}$ (which is also known as {\it
Kazarnovskii's pseudovolume \cite{kazarnovskii}}) is a multiple of the polynomial $f_n$ from Theorem
\ref{thm_fu}. This follows from the fact that the kernel of the restriction map
$\Val_n^{\mathrm{U}(n)} \to \Val_n^{\mathrm{U}(n-1)}$ is $1$-dimensional and contains
$\mu_{n,0}$ and $f_n$. 

With some more tricks, one can compute the scaling factor and compute the relations between the hermitian intrinsic volumes and the $ts$-basis. The result can be most easily expressed in terms of the Tasaki valuations: 
\begin{equation}
 \tau_{k,q} = \frac{\pi^k}{\omega_k(k-2q)! (2q)!} t^{k-2q}(4s-t^2)^q.
\end{equation}

Since $\Val^{\mathrm{U}(n)}$ is a finite-dimensional $\mathfrak{sl}_2$-representation, it admits a canonical decomposition ({\it Lefschetz decomposition}). The corresponding basis is called the {\it primitive basis} and is defined for all $0 \leq r \leq \frac {\min(k,2n-k)}{2}$ by   

\begin{equation} \label{eq_pi_in_tau}
\pi_{k,r}= (-1)^r(2n-4r+1)!!\sum_{i=0}^r (-1)^{i}\frac{(k-2i)!}{(2r-2i)!} \frac{(2r-2i-1)!!}{(2n-2r-2i+1)!!}\,\tau_{k,i}.
\end{equation}

This new basis is quite helpful for computational purposes, since in this basis, the matrix describing the Alesker-Poincar\'e-duality is anti-diagonal and we can easily compute its inverse (which is what we have to do in order to compute $k_{\mathrm{U}(n)}(\chi)$, see Theorem \ref{thm_ftaig}). 

As a result, the principal kinematic formula $k_{\mathrm{U}(n)}(\chi)$ in terms of the primitive basis was established in \cite{bernig_fu_hig}.

\begin{Theorem}\label{pkf}
Set $p:=\min\left\{\lfloor \frac k 2\rfloor ,\lfloor \frac{2n-k}2\rfloor\right\}$. 
\begin{align}\label{pkf1} 
 k_{U(n)}(\chi) &=\frac 1{\pi^n}\sum_{k=0}^{2n}{\omega_k\omega_{2n-k}}\\
\notag &\quad\quad\sum_{r=0}^p\frac{(n-r)!}{8^r(2n-4r)!} \frac{(2n-2r+1)!!}{(2n-4r+1)!!}\binom n {2r}^{-1}\pi_{k,r} \otimes{\pi_{2n-k,r}}
\end{align}
\end{Theorem}

Using \eqref{eq_pi_in_tau} and \eqref{eq_relation_tasaki_hiv}, we may restate this formula in terms of the Tasaki basis or in terms of intrinsic volumes, but the corresponding formulas are rather difficult. 

In conclusion, the vector space structure as well as the algebra structure on $\Val^{\mathrm{U}(n)}$ in terms of the different bases are now rather well understood. 

\subsection{Positive and monotone cone}

A valuation is called {\it positive} if $\mu(K) \geq 0$ for all $K$. It is easy to see that an $\mathrm{SO}(n)$-invariant valuation $\sum c_k\mu_k$ is positive if and only if all $c_k$ are positive. In fact, $\mu$ evaluated at a $k$-dimensional disk of radius $r$ behaves like $c_kr^k+o(r^k)$. On the other hand, it clearly follows from \eqref{eq_int_vol_projections} or from \eqref{eq_int_vol_intersections} that each $\mu_k$ is positive. Moreover, the $\mu_k$ and each positive linear combination $\mu$ of them is {\it monotone}, i.e. $\mu(K) \leq \mu(L)$ if $K \subset L$. Hence in the classical setting, the cones of positive and monotone invariant valuations coincide.

The situation in the $\mathrm{U}(n)$-case is more involved. A similar argument as above shows that a valuation $\mu=\sum_{k,q} c_{k,q}\mu_{k,q}$ can only be positive if the Klain function of each homogeneous component is positive, hence $c_{k,q} \geq 0$. That the $\mu_{k,q}$ are indeed positive does not follow immediately from their definition. But it can be shown (using the fact that the $\mu_{k,q}$ are {\it constant coefficient valuations}) that $\mu_{k,q}$ evaluated at a polytope is positive from which the positivity of $\mu_{k,q}$ follows by continuity.

What about the monotone cone? One way to construct an invariant monotone valuation is to use a positive invariant Crofton measure. It is not hard to see that the cone of all invariant valuations admitting a positive Crofton measure is dual to the positive cone with respect to the scalar product $\langle \phi,\psi\rangle:=\mathrm{PD}(\phi)(\hat \psi)$. 

But there are more monotone valuations. The idea to test monotonicity of a smooth, translation invariant valuation $\mu$ is to use a variation of a smooth convex body and to describe the first variation $\delta \mu$ of $\mu$ as a {\it curvature measure}, which is a signed measure concentrated on the boundary of $K$. 

The main observation is that $\mu$ is monotone if and only if the corresponding curvature measure is positive, and that this happens if and only if some {\it infinitesimal valuations} associated to $\delta \mu$ are positive. 

In the $\mathrm{U}(n)$-case, Park \cite{park02} has written down a list of equivariant curvature measures. The first variation map $\delta$ may be computed in terms of the hermitian intrinsic volumes and Park's curvature measures.  

Since the positive cone in $\Val^{\mathrm{U}(n)}$ is known (see above), we can thus determine the monotone cone too. The result is that a valuation 
\begin{displaymath}
\mu=\sum_{k,q} c_{k,q}\mu_{k,q}
\end{displaymath}
is monotone if and only if 
\begin{displaymath} 
(k-2q)c_{k,q} \geq (k-2q-1)c_{k,q+1}, \quad \max\{0,k-n\} \leq q \leq \left\lfloor \frac{k-1}{2}\right\rfloor
\end{displaymath}
and 
\begin{displaymath} 
(n+q-k+1) c_{k,q} \leq (n+q-k+3/2) c_{k,q+1}, \quad \max\{0,k-n-1\} \leq q \leq \left\lfloor \frac{k-2}{2}\right\rfloor.
\end{displaymath}

From this description we see that $\mu \in \Val^{\mathrm{U}(n)}$ is monotone if and only if each homogeneous component of $\mu$ is monotone. This is a general fact \cite{bernig_fu_hig}: {\it A translation invariant continuous valuation is monotone if and only if each homogeneous component is monotone}. This answers a question of P. McMullen \cite{mcmullen77}. The corresponding statement with {\it monotone} replaced by {\it positive} seems to be unknown. 

We can draw some more consequences of the above result. The cone of monotone invariant valuations is a polyhedral cone. It is not closed under any of the algebraic constructions from Section \ref{sec_alg_structures}. 
Let us give some examples (which are extremal rays of the polyhedral cone of monotone invariant valuations). 

The valuation 
\begin{displaymath}
 \mu:=\mu_{4,1}+\frac{2}{3}\mu_{4,2} \in \Val_4^{\mathrm{U}(3)}
\end{displaymath}
is monotone, but its Fourier transform 
\begin{displaymath}
 \hat \mu=\mu_{2,0}+\frac{2}{3}\mu_{2,1} \in \Val_2^{\mathrm{U}(3)}
\end{displaymath}
is not monotone (the second inequality with $q=0$ is violated). 

Consider 
\begin{align*}
 \mu & :=\mu_{4,0}+\frac{6}{7}\mu_{4,1}+\frac{12}{7}\mu_{4,2} \in \Val_4^{\mathrm{U}(6)}\\
 \phi & :=\mu_{4,0}+\frac{4}{3}\mu_{4,1}+\frac{32}{27}\mu_{4,2} \in \Val_4^{\mathrm{U}(6)}.
\end{align*}
Then $\mu,\phi$ are monotone valuations. From the technique described in Subsection \ref{subsec_kun}, one obtains that 
\begin{displaymath}
 \mu \cdot \phi=\frac{1002}{81}\mu_{8,2}+\frac{2552}{189}\mu_{8,3}+\frac{6112}{567}\mu_{8,4} \in \Val_8^{\mathrm{U}(6)},
\end{displaymath}
which is not monotone (the second inequality with $q=3$ is violated). 

Similarly, the invariant valuations 
\begin{align*}
 \mu & :=\mu_{4,0}+\frac{2}{3}\mu_{4,1}+\frac{4}{3}\mu_{4,2} \in \Val_4^{\mathrm{U}(4)}\\
 \phi & :=\mu_{6,2}+\frac{2}{3}\mu_{6,3} \in \Val_6^{\mathrm{U}(4)}.
\end{align*}
are monotone, but their convolution product
\begin{displaymath}
 \mu * \phi=4\mu_{2,0}+\frac{8}{3}\mu_{2,1} \in \Val_2^{\mathrm{U}(4)},
\end{displaymath}
is not monotone (the second inequality with $q=0$ is violated). 

This can be used to show that a monotone version of McMullen's conjecture does not hold true. Taking linear combinations of valuations of the form $K \mapsto \vol(K+A)$ with {\it positive} coefficients clearly yields monotone valuations and one would expect that every monotone valuation is the limit of such positive linear combinations. But this would imply that the monotone cone is closed under convolution, which is not the case.
 
\section{Other group actions}
\label{sec_other_groups}

Let $G$ be any compact connected Lie group acting transitively on the unit sphere. We have seen that $\Val^G$ is a finite-dimensional algebra, and that there are kinematic and additive $G$-kinematic formulas. Groups with this property are listed in \eqref{eq_transitive_groups} and \eqref{eq_exceptional_groups}. The classical case $G=\mathrm{SO}(n)$ was sketched in Section \ref{sec_classical_intgeo}, while the case $G=\mathrm{U}(n)$ was the subject of Section \ref{sec_un}. 

In this section, we will explain what is known for other $G$. 

\subsection{Special unitary group} 

The difference between the integral geometry of $\mathrm{U}(n)$ and that of $\mathrm{SU}(n)$ is not large. Naturally enough, it comes from the complex determinant. 

Let $V$ be a hermitian vector space of dimension $n$, and let $\mathrm{SU}(V) \cong \mathrm{SU}(n)$ be the special unitary group acting on $V$. 

For $k \neq n$, two $k$-dimensional subspaces are in the same $\mathrm{SU}(n)$-orbit if and only if they are in the same $\mathrm{U}(n)$-orbit. Klain's theorem thus implies that if $\mu \in \Val_k^+(V)$ is even and $\mathrm{SU}(n)$-invariant, then it is already $\mathrm{U}(n)$-invariant. 

As it turns out, all $\mathrm{SU}(n)$-invariant valuations are even. This is not trivial if $n \equiv 1 \mod 2$, since in this case $-1 \not\in \mathrm{SU}(n)$. 

In the middle degree however, things are different. Given an $n$-dimensional subspace $W$ in a complex $n$-dimensional vector space, one defines 
\begin{displaymath}
 \Theta(W):=\det(w_1,\ldots,w_n),
\end{displaymath}
where $w_1,\ldots,w_n$ is an orthonormal basis of $W$. Since another choice of basis $w_1,\ldots,w_n$ will affect $\Theta$ by the factor $\pm 1$ (depending on the orientations), this invariant is a well-defined element of $\mathbb{C}/\{\pm 1\}$. If the restriction of the symplectic form of $V$ on $W$ is not degenerated (which can only happen if $n$ is even), there is a natural choice of orientation of $W$ and $\Theta(W)$ is well-defined in $\mathbb{C}$. 

Two $\mathrm{U}(n)$-equivalent $n$-dimensional subspaces in $V$ belong to the same $\mathrm{SU}(n)$-orbit if and only if their $\Theta$-invariants agree. Using this, one can show that  
\begin{align*}
 \dim \Val_k^{\mathrm{SU}(n)}=\left\{\begin{array}{l l} \dim \Val_k^{\mathrm{U}(n)} & k \neq n\\ \dim \Val_k^{\mathrm{U}(n)}+4 & k=n, n \equiv 0 \mod 2\\ \dim \Val_k^{\mathrm{U}(n)}+2 & k=n, n \equiv 1 \mod 2.\end{array}\right.
\end{align*}
The Klain functions of the new valuations may be explicitly described in terms of Tasaki angles and the $\Theta$-invariant. The algebra structure and the kinematic formulas for $\mathrm{SU}(n)$ are variations from the $\mathrm{U}(n)$-case, see \cite{bernig_sun09}.

\subsection{Exceptional groups}

The group $\mathrm{Spin}(9)$ is the universal (two-fold) cover of $\mathrm{SO}(9)$. It can be explicitly described in a number of ways, for instance using Clifford algebras or using octonions. It acts on a $16$-dimensional space $\mathbb{R}^{16}$ which may be interpreted as an octonionic plane $\mathbb{O}^2$. The group $\mathrm{Spin}(7)$ acts on $\mathbb{R}^8$, which is an octonionic line. The group of automorphisms of $\mathbb{O}$ is called $\mathrm{G}_2$, it acts on the space of purely octonionic elements, which is $\mathbb{R}^7$.

Let $v$ be a point of the corresponding unit sphere. The stabilizers of $\mathrm{Spin}(9), \mathrm{Spin}(7)$ and $\mathrm{G}_2$ are given by $\mathrm{Spin}(7), \mathrm{G}_2$ and $S\mathrm{U}(3)$. The action of $\mathrm{G}_2$ on $T_vS^7$ and that of $\mathrm{SU}(3)$ on $T_vS^6$ are again transitive on the corresponding unit spheres, which is not the case for the action of $\mathrm{Spin}(7)$ on $T_vS^{15}$. This makes it rather easy to describe the integral geometry of $\mathrm{G}_2$ and $\mathrm{Spin}(7)$, but for $\mathrm{Spin}(9)$ other methods will be necessary. 

Let us first consider $\mathrm{G}_2$. The stabilizer is $\mathrm{SU}(3)$ acting on $W:=T_vS^6$. Any $\mathrm{G}_2$-invariant valuation $\mu$ may be restricted to a $\mathrm{SU}(3)$-invariant valuation on $W$. The restriction of $\mu$ to $W$ vanishes if and only if $\mu$ is simple, since $\mathrm{G}_2$ acts transitively on $6$-dimensional subspaces. But simple valuations are of degree $7$ (in the even case) or $6$ (in the odd case). Hence, if $\mu$ is of degree $k \leq 5$, $\mu|_W=0$ if and only if $\mu=0$, and therefore $\dim \Val_k^{\mathrm{G}_2} \leq \dim \Val_k^{\mathrm{SU}(3)}$ for $0 \leq k \leq 5$. Using furthermore the symmetry induced by the Hard Lefschetz theorem, we obtain that $\dim \Val_k^{\mathrm{G}_2}=1$ for $k \neq 3,4$ and that $\dim \Val_3^{\mathrm{G}_2}=\dim \Val_4^{\mathrm{G}_2}$ is either $1$ or $2$. 

Now we repeat the argument with $\mathrm{Spin}(7)$ instead of $\mathrm{G}_2$ and $\mathrm{G}_2$ instead of $\mathrm{SU}(3)$ to obtain that $\dim \Val_k^{\mathrm{Spin}(7)}=1$ for $k \neq 4$ and that $\dim \Val_4^{\mathrm{Spin}(7)}$ equals $1$ or $2$. 

It remains to decide whether $\dim \Val_4^{\mathrm{Spin}(7)}$ equals $1$ or $2$. Since $\mathrm{Spin}(7)$ contains $\mathrm{SU}(4)$ as a subgroup, it is easy to find a $\mathrm{Spin}(7)$-invariant, not $\mathrm{SO}(8)$-invariant element of degree $4$ in $\Val^{\mathrm{SU}(4)}$. Going back, we see that this implies $\dim \Val_3^{\mathrm{G}_2}=\dim \Val_4^{\mathrm{G}_2}=2$, hence we get the following table: 
\begin{displaymath}
\begin{array}{c | c | c} k & \dim \Val_k^{\mathrm{G}_2} & \dim \Val_k^{\mathrm{Spin}(7)} \\ \hline
0 & 1 & 1\\
  1 & 1 & 1\\
2 & 1 & 1 \\
3 & 2 & 1\\
4 & 2 & 2 \\
5 & 1 & 1\\
6 & 1 & 1\\
7 & 1 & 1\\
8 & - & 1
\end{array}
\end{displaymath}

The new valuations in these spaces may be explicitly described. Since there are relatively few of them valuations, it is an easy task to describe the product structures. There are isomorphisms of graded algebras 
\begin{align*}
 \Val^{\mathrm{G}_2} & \cong \mathbb{C}[t,u]/(t^2u,u^2+4t^6)\\
\Val^{\mathrm{Spin}(7)} & \cong \mathbb{C}[t,v]/(v^2-t^8,vt),
\end{align*}
where $u$ is of degree $3$ and $v$ of degree $4$. 

From these isomorphisms and the fundamental theorem of algebraic integral geometry, one can derive kinematic formulas and additive formulas for $\mathrm{G}_2$ and $\mathrm{Spin}(7)$. We refer to \cite{bernig_g2} for details.

\subsection{Symplectic groups}
\label{sub_sec_sympgroups}

The integral geometry of the remaining three sequences $\mathrm{Sp}(n), \mathrm{Sp}(n) \cdot \mathrm{U}(1), \mathrm{Sp}(n) \cdot \mathrm{Sp}(1)$ in the list \eqref{eq_transitive_groups} seems to be quite difficult. The case $n=1$ is already contained in the $\mathrm{SU}(n)$-theory, since $\mathrm{Sp}(1) \cong \mathrm{SU}(2)$ (see also \cite{alesker_su2_04,bernig_quat09}). But even for $n=2$, things are mysterious. From a combinatorial formula in \cite{bernig_qig}, one gets the following dimensions

\begin{displaymath}
 \begin{array}{c | c c c c c c c c c }  k & 0 & 1 & 2 & 3 & 4 & 5 & 6 & 7 & 8\\ \hline
\dim \Val_k^{\mathrm{Sp}(2)} & 1 &  1 & 7 & 13 & 29 & 13 & 7 & 1 & 1\\
\dim \Val_k^{\mathrm{Sp}(2)\cdot \mathrm{U}(1)} & 1 & 1 & 3 & 5 & 9 & 5 & 3 & 1 & 1\\
\dim \Val_k^{\mathrm{Sp}(2) \cdot \mathrm{Sp}(1)} & 1 & 1 & 2 & 3 & 5 & 3 & 2 & 1 & 1
 \end{array}
\end{displaymath}

It is not known how to describe these valuations geometrically. For general $n$, there is a combinatorial formula using Young diagrams and Schur functions to compute the dimensions of the spaces $\Val_k^{\mathrm{Sp}(n)}, \Val_k^{\mathrm{Sp}(n)\cdot \mathrm{U}(1)}, \Val_k^{\mathrm{Sp}(n) \cdot \mathrm{Sp}(1)}$. The behavior of these numbers is rather irregular. For large $n$ (in fact $n \geq k$ is enough), these dimensions stabilize to some value $\dim \Val_k^{\mathrm{Sp}(\infty)}$ (resp. $\dim \Val_k^{\mathrm{Sp}(\infty) \cdot \mathrm{U}(1)}$, $\dim \Val_k^{\mathrm{Sp}(\infty) \cdot \mathrm{Sp}(1)}$). These asymptotic values can be explicitly computed, their Poincar\'e series is given by 

\begin{align*} 
\sum_{k=0}^\infty \dim \Val_k^{\mathrm{Sp}(\infty)} x^k &
=\frac{x^4-3x^3+6x^2-3x+1}{(1-x)^7(1+x)^3}\\
\sum_{k=0}^\infty \dim \Val_k^{\mathrm{Sp}(\infty)\cdot \mathrm{U}(1)} x^k &
=\frac{x^6-2x^5+2x^4+2x^2-2x+1}{(x^2+1)(x^2+x+1)(1+x)^2(1-x)^6}\\
 \sum_{k=0}^\infty \dim \Val_k^{\mathrm{Sp}(\infty)\cdot \mathrm{Sp}(1)} x^k &
=\frac{x^5+2x^4+x^3+1}{(x^2+1)(x^2+x+1)(1+x)^2(1-x)^4}.
\end{align*}

This is another hint that quaternionic integral geometry is difficult, since it follows from these expressions that none of the algebras $\Val^{\mathrm{Sp}(\infty)}, \Val^{\mathrm{Sp}(\infty) \cdot \mathrm{U}(1)}, \Val^{\mathrm{Sp}(\infty) \cdot \mathrm{Sp}(1)}$ is a freely generated algebra (in contrast to  the $\mathrm{U}(n)$-case, where $\Val^{\mathrm{U}(\infty)} \cong \mathbb{C}[t,s]$). It is not even known whether these algebras are finitely generated.

\section{Some open problems}
\label{sec_open_problems}

Let us describe three main problems whose solutions will probably stimulate further progress in algebraic integral geometry. 

\begin{enumerate}
 \item In the Euclidean setting, there are more elaborate versions of the kinematic formulas, the {\it local} kinematic formulas \cite{schneider_book93}. They apply to {\it curvature measures}, which are local versions of the intrinsic volumes. Each intrinsic volume is related to exactly one curvature measure. In the hermitian case, the invariant curvature measures were described in \cite{bernig_fu_hig}. It is known that there are local kinematic formulas \cite{fu90}. However, the computation of the coefficients in such a formula is a challenge, since the algebraic machinery from Section \ref{sec_alg_structures} only applies to valuations and not to curvature measures. 
\item We have described in detail the theory of valuations on an affine space. The theory of {\it valuations on manifolds} was recently worked out, mainly by Alesker \cite{alesker_val_man1, alesker_val_man2, alesker_val_man3, alesker_val_man4, alesker_survey07, bernig_broecker07, bernig_quat09, alesker_bernig}. This gives the appropriate framework to study integral geometry of projective and hyperbolic spaces. It turns out that on compact rank one symmetric spaces (CROSS), the space of (smooth) invariant valuations is finite-dimensional and that a version of the fundamental theorem of algebraic integral geometry holds true \cite{alesker_bernig}. To work out the algebraic structure of the space of valuations on a CROSS is a challenge. In the case of $\mathbb{CP}^n$, Abardia \cite{abardia_thesis} and Abardia-Gallego-Solanes \cite{abardia_gallego_solanes} studied various Crofton- and Chern-Gauss-Bonnet-type formulas. There also exists a (rather mysterious) conjecture by J. Fu \cite{oberwolfach_report} concerning the algebra structure of the space of invariant valuations on $\mathbb{CP}^n$.
\item The intrinsic volumes satisfy a number of important inequalities, like the isoperimetric inequality and the Brunn-Minkowski inequality. What is the corresponding statement in the hermitian case? A special case of this general question is the following. Let $\mu:=a\mu_{2,0}+b\mu_{2,1} \in \Val_2^{\mathrm{U}(2)}$ be a positive valuation. What is the minimum of $\mu(K)$ as $K$ ranges over all compact convex bodies of volumes $1$? By a version of the isoperimetric inequality, the minimum in the case $a=b$ is achieved by a ball, but the case $a \neq b$ is open.         
\end{enumerate}


\end{document}